\documentclass[12pt]{article}

\usepackage{amsmath}
\usepackage{amssymb}
\usepackage{latexsym}
\usepackage[raiselinks=false,colorlinks=true,citecolor=blue,urlcolor=blue,linkcolor=blue,bookmarksopen=true,dvips]{hyperref}
\usepackage{diagrams}
\diagramstyle{h=21pt,w=21pt,scriptlabels,midshaft,PostScript=dvips,nohug,shortfall=4pt}

\newcommand{\name}[1]{\ulcorner #1 \urcorner}

\usepackage[T1]{fontenc}

\usepackage[noBBpl]{mathpazo}
\usepackage{graphicx}
\renewcommand{\texttt}[1]{{\fontfamily{pcr}\fontseries{m}\fontshape{n}%
\selectfont #1}}

\usepackage{mathrsfs}
\providecommand{\lastUpdate}[1]{#1}

\newcommand{\FinSet}{\kat{FinSet}}

\newcommand{\tr}{\mathbf{T}}

\newcommand{\mybialg}{\mathscr{B}}
\newcommand{\CK}{\mathscr{H}}

\newcommand{\isopil}{\stackrel{\raisebox{0.1ex}[0ex][0ex]{\(\sim\)}}%
			{\raisebox{-0.15ex}[0.28ex]{\(\rightarrow\)}}}

\newcommand{\upperstar}{^{\raisebox{-0.25ex}[0ex][0ex]{\(\ast\)}}}
\newcommand{\lowerstar}{_{\raisebox{-0.33ex}[-0.5ex][0ex]{\(\ast\)}}}
\newcommand{\lowershriek}{_!}
\newcommand{\upperdot}{^\bullet}
\newcommand{\df}{\: {\raisebox{0.255ex}{\normalfont\scriptsize :\!\!}}=}

\newcommand{\tensor}	{\otimes}
\newcommand{\ground}{\Bbbk}
\newcommand{\Spec}{\operatorname{Spec}}
\newcommand{\Id}{\operatorname{Id}}
\newcommand{\Aut}{\operatorname{Aut}}

\newcommand{\isleftadjointto}{\dashv}
\newcommand{\op}{^{\text{{\rm{op}}}}}

\newcommand{\DD}{\mathscr{D}}

\providecommand{\kat}[1]{\text{\textbf{\textsl{#1}}}}

\newcommand{\Set}{\kat{Set}}

\newcommand{\A}{\mathbb{A}}
\newcommand{\B}{\mathbb{B}}
\newcommand{\F}{\mathbb{F}}
\newcommand{\N}{\mathbb{N}}
\newcommand{\Z}{\mathbb{Z}}

\newcommand{\ov}{\overline}
\newcommand{\fy}{\varphi}
\newcommand{\wtil}{\widetilde}

\newtheorem{taller}[subsection]{$\!\!$}
\newenvironment{blanko}[1]%
{\begin{taller}{\normalfont\bfseries #1}\normalfont}%
{\end{taller}}

\newcommand{\overskrift}[1]{\par\noindent\relax{\LARGE #1}\par\bigskip}

\usepackage{texdraw}
\input{txdtools}

\everytexdraw{%
	\drawdim pt \linewd 0.5 \textref h:C v:C
	\setunitscale 1
}

\newcommand{\onedot}{
  \bsegment
    \move (0 0) \fcir f:0 r:2
  \esegment
}

  \newcommand{\inlineDotlessTree}{%
\raisebox{-4pt}{
\begin{texdraw} \linewd 0.5 \bsegment
    \move (0 0) \lvec (0 15) \move (1 0)
  \esegment \end{texdraw} } }

  \newcommand{\smalltree}{
  \bsegment
  \linewd 1.5
    \move (0 0) \fcir f:0 r:3
    \lvec (-10 20) \fcir f:0 r:3
    \move (0 0) \lvec (10 20) \fcir f:0 r:3
    \lvec (0 40) \fcir f:0 r:3
    \move (10 20) \lvec (20 40) \fcir f:0 r:3
    \savepos (0 0)(*ex *ey)
  \esegment
  \move (*ex *ey)
}

\newcommand{\allleaves}{
  \bsegment
  \linewd 0.5
    \move (0 40) \lvec (-5 65) 
    \move (0 40) \lvec (5 65) 
    \move (0 0) \lvec (0 -15)
    \move (0 0) \lvec (-40 45)
    \move (-10 20) \lvec (-22 60)
    \move (10 20) \lvec (15 64)
    \savepos (0 0)(*ex *ey)
  \esegment
  \move (*ex *ey)
}

\newcommand{\Poly}{\kat{Poly}}

\newarrow{W}{=}{=}{=}{=}{>>}
\newarrow{Iso}{>}{=}{=}{=}{>>}
\newarrow{Lig}{=}{=}{=}{=}{=}
\newarrow{Mono}{>}{-}{-}{-}{>}
\newarrow{Epi}{-}{-}{-}{-}{>>}
\newarrow{Face}{>}{-}{-}{-}{>}
\newarrow{Deg}{-}{-}{-}{-}{>>}
\newcommand{\Aff}{\kat{Aff}}
\newcommand{\TEmb}{\kat{TEmb}}

\makeatletter
\renewcommand{\ps@headings}
	{\setlength{\headheight}{13pt}%
	 \setlength{\headsep}{12pt}%
	 \renewcommand{\@oddhead}{\parbox{\textwidth}{%
			\footnotesize
			\texttt{\jobname.tex \ \ \lastUpdate{2012-04-28 14:57} \ \ \ 
			\hfill [\thepage/\pageref{lastpage}]}
			\\ \rule[8pt]{\textwidth}{0.3pt}}%
	 }
	\renewcommand{\@oddfoot}{}
	\renewcommand{\@evenfoot}{}%
}
\makeatother

\usepackage{color}

\newcommand{\putsection}{
\begin{picture}(8,1)(0,0)	
  \put(0,0){\line(1,0){8}}
\end{picture}
}
\newcommand{\putshortsection}{
\begin{picture}(8,1)(0,0)
  \put(0,0){\line(1,0){3.5}}
\end{picture}
}

\newcommand{\putpointedsection}{
\begin{picture}(8,2)(0,0)	
  \put(2,-1.7){\line(0,1){3.8}}
  \put(0,0){\line(1,0){8}}
\end{picture}
}

\newcommand{\putsectionsection}{
\begin{picture}(8,0)(0,0)	
  \put(0,0){\line(1,0){8}}
  \put(0,2){\line(1,0){8}}
\end{picture}
}

\newcommand{\putsectionpointedsection}{
\begin{picture}(8,0)(0,0)
  \color{blue}
  \put(0,0){\line(1,0){8}}
  \put(2,-1.7){\line(0,1){3.8}}
  \color{red}
  \put(0,2){\line(1,0){4}}
  \color{black}
\end{picture}
}

\newcommand{\putredsectionbluesection}{
\begin{picture}(8,0)(0,0)
  \color{blue}
  \put(0,0){\line(1,0){8}}
  \color{red}
  \put(0,2){\line(1,0){8}}
  \color{black}
\end{picture}
}

\newcommand{\sect}[1]{\overset{\putsection}{#1}}
\newcommand{\ptsect}[1]{\overset{\putpointedsection}{#1}}
\newcommand{\sectsect}[1]{\overset{\putsectionsection}{#1}}
\newcommand{\sectptsect}[1]{\overset{\putsectionpointedsection}{#1}}

\newcommand{\shortsect}[1]{\overset{\putshortsection}{#1}}

\newcommand{\red}[1]{\color{red}{#1}\color{black}}
\newcommand{\blue}[1]{\color{blue}{#1}\color{black}}
\newcommand{\redsect}[1]{\overset{\red{\putsection}}{#1}}
\newcommand{\redptsect}[1]{\overset{\red{\putpointedsection}}{#1}}
\newcommand{\redsectbluesect}[1]{\overset{\putredsectionbluesection}{#1}}
\newcommand{\bluesect}[1]{\overset{\blue{\putsection}}{#1}}
\newcommand{\blueptsect}[1]{\overset{\blue{\putpointedsection}}{#1}}

\newcommand{\arxiv}[1]{\href{http://arxiv.org/pdf/#1}{ArXiv:#1}}

\begin{document}

\pagestyle{headings}

\newcommand{\polyFunct}[8]{
\begin{diagram}[w=2.5ex,h=3.7ex,tight]
&&#2&&\rTo^{#6}&&#3\\
&\ldTo^{#5}&&&&&&\rdTo^{#7}\\
#1&&&&#8&&&&#4
\end{diagram}
}

\vspace*{24pt}

\begin{center}
  
\overskrift{Categorification of Hopf algebras of rooted trees}

\bigskip

\textsc{Joachim Kock}\footnote{
Departament de Matem\`atiques,
Universitat Aut\`onoma de
Barcelona,
08193 Bellaterra,
Spain}

\end{center}

\begin{abstract}
  We exhibit a monoidal structure on the category of finite sets indexed by
  $P$-trees for a finitary polynomial endofunctor $P$.  This structure
  categorifies the monoid scheme (over $\Spec \N$) whose semiring of functions
  is (a $P$-version of) the Connes--Kreimer bialgebra $H$ of rooted trees
  (a Hopf algebra after base
  change to $\Z$ and collapsing $H_0$).  The monoidal structure is itself 
  given by a polynomial functor,
  represented by three easily described set maps; we show that these
  maps are the same as those occurring in the
  polynomial representation of the free monad on $P$.
\end{abstract}

\small



% 18A99 General theory of categories and functors -- None of the above, but in this section
% 05C05    	Trees
% 16T05    	Hopf algebras and their applications [See also 16S40, 57T05]
% 16B50    	Category-theoretic methods and results (except as in 16D90) [See also 18-XX]
% 18D10    	Monoidal categories (= multiplicative categories), symmetric monoidal categories, braided categories [See also 19D23]
% 	

\normalsize

%%%%%%%%%%%%%%%%%%%%%%%%%%%%%%%%%%%%%%%%%%%%%%%%%%
\section{Introduction}
%%%%%%%%%%%%%%%%%%%%%%%%%%%%%%%%%%%%%%%%%%%%%%%%%%

\begin{blanko}{The Connes--Kreimer Hopf algebra and the Butcher group.}
  The Hopf algebra $\CK$ of rooted trees is now a well-established object in
  mathematics, thanks in particular to the seminal works of Connes and Kreimer.
  Kreimer~\cite{Kreimer:9707029} discovered
  that $\CK$ controls the combinatorics of renormalisation in perturbative quantum
  field theory, and his collaboration with Connes 
  (e.g.~\cite{Connes-Kreimer:9808042}, \cite{Connes-Kreimer:9912092}, to cite a
  few) uncovered deep connections with noncommutative geometry, number
  theory, Lie theory, and algebraic combinatorics,  stimulating a lot of further
  activity by many mathematicians and physicists.  The Connes--Kreimer Hopf algebra has now
  been characterised by several different universal properties
  \cite{Connes-Kreimer:9808042}, \cite{Moerdjik:9907010},
  \cite{Chapoton-Livernet:0002069}.
  The group of characters of $\CK$, now called the Butcher
  group, was in fact studied by Butcher~\cite{Butcher:1972} some 30 years 
  earlier,
  in relation with order conditions for Runge--Kutta methods in numerical
  integration. 
%   D\"ur~\cite{Dur:1986} had realised in 1986 that the Hopf algebra
%   of representable functions of the Butcher group is in fact an incidence Hopf
%   algebra, and
  The link back to this work from the Connes--Kreimer Hopf algebra was
  provided by Brouder~\cite{Brouder:9904014}.  
  
%   It is now an important object, and it has been characterised by several
%   universal properties explaining its importance.  In quantum field theory, the
%   most important characterisation seems to be that it contains the initial
%   \ldots cocycle (the operator that adjoins a root to a forest), in terms of
%   Cartier--Hochschild cohomology.  Moerdijk~\cite{Moerdjik:9907010} showed that
%   it is the initial $\operatorname{Comm}[t]$-algebra, where
%   $\operatorname{Comm}[t]$ is the operad for commutative algebras, with a freely
%   adjoint unary operation $t$.  In combinatorics, the most important
%   characterisation is perhaps as the dual of the enveloping algebra of the free
%   pre-Lie algebra on one generator 
%   (Chapoton--Livernet~\cite{Chapoton-Livernet:0002069}).
  
% Schmitt incidence algebra of poset \cite{Schmitt:incidence}
% See also Loday--Ronco~\cite{Loday-Ronco:0810.0435}.
\end{blanko}

\begin{blanko}{Categorification of the Connes--Kreimer Hopf algebra.}
  This article investigates a categorification of the bialgebra $\CK$, i.e.~a lift
  from the level of algebra to the level of sets.  To be precise, the antipode
  is {\em not} categorified.  Indeed, the antipode is a feature depending on the
  additive inverses of the ground ring.  But in fact, the most natural `ring' of
  definition of $\CK$ is the semiring $\N$ of natural numbers, and lacking
  additive inverses there is no
  antipode to be had here.  The semiring $\N$ appears as the Burnside semiring of the
  distributive category $\F$ of finite sets, i.e.~is the set of isomorphism
  classes of finite sets, with addition and multiplication inherited from the
  categorical sum and product.  Exhibiting a distributive category whose
  Burnside semiring is $\CK$, and describing at this level the comultiplication,
  is what we mean by categorification, following a popular terminology which
  goes back to Crane and Yetter~\cite{Crane-Yetter:categorification}; the
  specific process employed can more precisely be called objectification,
  cf.~Lawvere, Schanuel, and their collaborators.  A recommended introduction
  to categorification is the beautiful paper \cite{Baez-Dolan:finset-feynman} of
  Baez and Dolan.
  
  While the combinatorial nature of $\CK$ makes it clear that this
  categorification should be possible, a considerable amount of categorical
  algebra is needed in order to make all the algebraic structure explicit at the
  objective level.  The categorification of $\CK$ will be the distributive category
  $\F[\tr]$ of polynomial functors on the category $\A^\tr$ of $\tr$-indexed
  finite sets, where $\tr$ is the set of trees.  For technical reasons we mostly
  work with $P$-trees for a finitary polynomial endofunctor $P$, and we start
  out by explaining the differences.  The comultiplication will be a 
  comomonoidal structure on $\F[\tr]$, relative to a certain tensor product 
  of distributive categories.  It is conceptually much easier to take the dual 
  viewpoint.
  The category $\A^\tr$ will be the
  categorification of the Butcher group (or rather the Butcher monoid), and the
  task is to describe the monoidal structure on $\A^\tr$ dual to the
  comultiplication.  This monoidal structure, a functor $M:\A^\tr \times \A^\tr
  \to \A^\tr$, is itself a polynomial functor.  The comultiplication is now
  given by precomposition with $M$.
\end{blanko}

\begin{blanko}{Polynomial functors.}    
  The notion of polynomial functor is central to this work.  The theory of
  polynomial functors has roots in topology, representation theory,
  combinatorics, logic and computer science, but the task of unifying these
  developments has only recently begun~\cite{Gambino-Kock:0906.4931}.  An
  important feature of polynomial functors is that they can be manipulated in
  terms of a few representing sets, just as polynomial functions can be
  manipulated in terms of their coefficients and exponents.
  
  In the
  present work, polynomial functors enter at two levels: firstly, and most
  importantly, the notion of polynomial functor categorifies the notion of
  polynomial function; second, the trees that index the involved variables
  are operadic trees, and they are themselves defined in terms of polynomial
  endofunctors.  For the first aspect we need to develop some basic theory,
  which constitutes Section~\ref{Sec:alg}; for the second, the theory needed is
  already available from \cite{Gambino-Kock:0906.4931} and \cite{Kock:0807}.
\end{blanko}

\begin{blanko}{Free monads, $P$-trees, and beyond.}
  Historically, one starting point for the general project of categorification
  is the quest in combinatorics for bijective proofs: it is well appreciated
  that a bijection between sets represents better understanding than a mere
  equation between numbers.
  One insight into the Hopf algebra of rooted trees which results from its
  categorification and the polynomial viewpoint is the relation with free
  monads: it is shown that the set maps occurring in the polynomial
  representation of the new monoidal structure are the same maps as occur in the
  polynomial representation of the free monad construction.  Unfortunately, this
  is not completely true for abstract trees.  It is true for $P$-trees, and we
  develop the theory in this setting --- the free monad in question is then the
  free monad on $P$.  The notion of $P$-tree (cf.~\ref{Ptree} below) covers many
  notions of structured and decorated trees, such as binary and planar
  trees, but abstract trees themselves are not an example: abstract trees should be
  $P$-trees for the terminal polynomial endofunctor, but the category of
  polynomial endofunctors over sets (see \ref{cart}) does not have a terminal object!

  Nevertheless, with a little care, the constructions made for $P$-trees work
  also for abstract trees, only the relation with free monads is then less
  direct.  This is explained in Section~\ref{Sec:beyond}, where it is also
  explained how these issues disappear when upgrading the whole theory from sets
  to groupoids: in this fancier setting, the terminal polynomial endofunctor
  does exist, and abstract trees become a particular case of the notion of
  $P$-tree.
\end{blanko}

\begin{blanko}{Outline of the paper.}
  In Section~\ref{Sec:bialg} we introduce the bialgebra of operadic trees, and
  explain its relation with the usual Connes--Kreimer Hopf algebra.  
  Section~\ref{Sec:poly} collects the notions and results needed from the theory
  of polynomial functors, notably the explicit formula for substitution of 
  polynomials.  In Section~\ref{Sec:alg} we set up a framework for dealing with
  categories of polynomial functors and categories of indexed finite sets as
  `polynomial rings' and `affine space', respectively.  This section contains
  many observations that seem not to have been made before, 
  regarding polynomial functors as categorification of polynomial functions.
  In Section~\ref{Sec:trees} we introduce trees and $P$-trees and review in
  detail the construction of the free monad on a polynomial endofunctor $P$.
  Finally in Section~\ref{Sec:monoidal} we establish the monoidal structure
  on $\A^\tr$, and relate it to the bialgebra of trees.  Section~\ref{Sec:beyond}
  contains some remarks about the difference between $P$-trees and abstract
  trees, and hints at a groupoid version of all the constructions as a more 
  comprehensive framework.
\end{blanko}

\begin{blanko}{Related work.}
  This work was presented at the 2010 International Category Theory Conference
  in Genova.  At the same conference, Mat\'\i as Menni spoke about his work with
  Lawvere~\cite{Lawvere-Menni} on categorification of incidence algebras of
  M\"obius categories.  While it is known that the Connes--Kreimer Hopf algebra
  can be constructed as the incidence algebra of a suitable family of posets, it
  seems unlikely that it can be given as the incidence algebra of a single
  M\"obius category, so at present the Lawvere--Menni approach does not apply 
  to categorify the Connes--Kreimer Hopf algebra.
  Work is under way to develop a higher-categorical notion of 
  M\"obius categories in order to
  unify the two approaches.
%   : in collaboration with Imma G\'alvez and Andy
%   Tonks~\cite{Galvez-Kock-Tonks:Moebius}, a theory of M\"obius monads is being
%   developed, and the bialgebra of rooted trees shown to be an example of an
%   incidence algebra of a M\"obius monad.  Briefly, just as every free category
%   is M\"obius, so is every free monad.
%   
%   In a related line of development, the explicit polynomial functor
%   constructions of the present paper are used in \cite{Galvez-Kock-Tonks:FdB} to
%   establish Fa\`a di Bruno formulae for Green functions \`a la van
%   Suijlekom~\cite{vanSuijlekom:0807}.
\end{blanko}

\begin{blanko}{Acknowledgments.}
  The author has benefited greatly from many conversations with Kurusch
  Ebrahimi-Fard, and also thanks the anonymous referees for many pertinent 
  remarks which led to improved exposition.
  Partial support from research grants
  MTM2009-10359 % Nart
  and MTM2010-20692 % Castellana
  of Spain and 2009SGR-1092 % Aguad\'e
  of Catalonia is gratefully acknowledged.
\end{blanko}

%%%%%%%%%%%%%%%%%%%%%%%%%%%%%%%%%%%%%%%%%%%%%%%%%%
\section{Hopf algebras of rooted trees}
%%%%%%%%%%%%%%%%%%%%%%%%%%%%%%%%%%%%%%%%%%%%%%%%%%

\label{Sec:bialg}

The standard Connes--Kreimer Hopf algebra of rooted trees concerns combinatorial
trees, whereas in this paper we prefer to work with operadic trees.
% , which have more expressive power.
This section explains the differences.

\begin{blanko}{Combinatorial trees.}
  The trees usually employed, which here we call combinatorial trees, are often
  defined as finite connected graphs without loops or cycles, and with a designated root
  vertex.  If a connected subgraph of a rooted tree $T$ does not contain the root of $T$, then 
  instead it has a vertex nearest the root, which is then defined to be its root.
\end{blanko}

\begin{blanko}{The Connes--Kreimer bialgebra of rooted trees.}
  The bialgebra of rooted trees of Connes and 
  Kreimer~\cite{Kreimer:9707029}   is the polynomial $\ground$-algebra $\CK$
  on the set of isomorphism classes of combinatorial trees.
  Here $\ground$ can be any commutative $\N$-algebra.
  The comultiplication is given on generators by
  \begin{eqnarray*}
  \Delta:  \CK & \longrightarrow & \CK \tensor_\ground \CK  \\
    T & \longmapsto & \sum_c P_c \tensor R_c ,
  \end{eqnarray*}
  where the sum is over all admissible cuts of $T$; the left-hand factor $P_c$
  is the forest (interpreted as a monomial) found above the cut, and $R_c$ is
  the subtree found below the cut (or the empty forest, in case the cut is below
  the root).  `Admissible cut' means upper-set in the poset underlying the tree
  (oriented from leaves (inputs) to root (output))
  i.e.~either a subtree containing the root, or the empty set.
%   (The notion of
%   cut is a source of much confusion in the literature, where for example some
%   authors define a cut to be a certain subset of the edges of the tree, and yet
%   claim to be able to distinguish the empty cut above the tree from the empty
%   cut below the tree.)

  $\CK$ is a connected bialgebra: the grading is by the number of nodes, and
  $\CK_0$ is spanned by the unit.  Therefore, by general principles (see for 
  example \cite{Figueroa-GraciaBondia:0408145}), it acquires
  an antipode and becomes a Hopf algebra --- provided $\ground$ has additive
  inverses, i.e.~is a $\Z$-algebra.  (In any case, the antipode exists after
  base change to $\Z$.)
\end{blanko}

\begin{blanko}{Operadic trees.}
  We shall need trees with slightly more expressive power, by allowing loose
  ends (leaves): these are {\em operadic trees}, also called finite rooted trees
  with boundary --- a formal definition is given in \ref{polytree-def}.  For the
  moment, the following drawings should suffice to exemplify operadic trees ---
  as usual the planar aspect inherent in a drawing should be disregarded:
  \begin{center}\begin{texdraw}
  \linewd 0.5 \footnotesize
  \move (-50 0)
  \bsegment
    \move (0 0) \lvec (0 30)
  \esegment
  
  \move (0 0)
  \bsegment
    \move (0 0) \lvec (0 18) \onedot
  \esegment
  
  \move (50 0)
  \bsegment
    \move (0 0) \lvec (0 36)
    \move (0 18) \onedot
  \esegment
  
  \move (105 0)
  \bsegment
    \move (0 0) \lvec (0 15) \onedot
    \lvec (-5 33) \onedot
    \move (0 15) \lvec (-12 28) \onedot
    \move (0 15) \lvec (4 43)
    \move (0 15) \lvec (12 40)
  \esegment
  
  \move ( 170 0)
  \bsegment
    \move (0 0) \lvec (0 18) \onedot
    \lvec (-6 32) \onedot
    \lvec (-12 57)
    \move (0 18) \lvec (4 40) \onedot
    \lvec (20 50) \onedot
    \lvec (15 65)
    \move (20 50) \lvec (25 65)
    \move (4 40) \lvec (9 54) \onedot
    \move (4 40) \lvec (-4 61)
  \esegment
  \end{texdraw}\end{center}
  Note that certain edges (the {\em leaves}) do not start in a node and that one
  edge (the obligatory {\em root edge}) does not end in a node.  A node without incoming
  edges is not the same thing as a leaf; it is a nullary operation (i.e.~a
  constant) in the sense of operads.  In operad theory, the nodes represent
  operations, and trees are formal combinations of operations.  The small
  incoming edges drawn at every node serve to keep track of the arities of the
  operations.  Furthermore, for coloured operads, the operations have type
  constraints on their inputs, encoded as attributes of the edges.
\end{blanko}

\begin{blanko}{Operadic trees in QFT.}
  The use in quantum field theory of more refined notions of trees, and operadic
  trees in particular, has been hinted at by Kreimer in several papers, most
  concretely with Bergbauer~\cite{Bergbauer-Kreimer:0506190}, where trees with
  loose edges are used to analyse combinatorial Dyson--Schwinger equations.

  In fact, the role of trees in the Connes--Kreimer bialgebra is to encode
  nestings of 1-particle irreducible Feynman graphs~\cite{Kreimer:9707029},
  and one can argue \cite{Kock:graphs-and-trees} that they are naturally
  operadic: each tree naturally comes equipped with
  decorations by primitive 1PI graphs on nodes and interaction labels on edges.
  To fully encode the compatibility conditions involved in these decorations,
  and to allow to recover the graph from the decorated trees, it is necessary to
  keep track of the arities of the nodes, so that even vacant input slots are
  represented; this leads naturally to operadic trees.  A more thorough analysis
  of the relationship between graphs and trees is given
  elsewhere~\cite{Kock:graphs-and-trees}.

  Kreimer~\cite{Kreimer:9810022} stresses the general utility of trees as
  expression of nestings of structures, not only of Feynman graphs.  In the
  following picture we see first a nesting of subsets, then a combinatorial-tree
  expression of the nesting, and finally an operadic-tree expression of the same
  nesting, in which the leaves correspond to the elements of the nested sets:
\begin{center}
\begin{texdraw}
  \linewd 0.5 \footnotesize
  \bsegment
    \move (0 30) \lcir r:30
    \move (-15 15) \lcir r:7
    \move (8 32) \lcir r:19
    \move (17 38) \lcir r:5 %\onedot
    \move (7 24) \lcir r:9
    \move (-15 15) \onedot
    \move (4 25) \onedot
    \move (10 23) \onedot
    \move (4 42) \onedot
%     \move (-19 36) \onedot
    \move (-13 48) \onedot
  \esegment
\end{texdraw}
\hspace{4em}
\begin{texdraw}
  \setunitscale 0.8
  \move (0 15)\smalltree
  \move (0 0)
\end{texdraw}
\hspace{3em}
\begin{texdraw}
  \setunitscale 0.8
  \smalltree\allleaves
%   \linewd 0.5 \footnotesize
%   \bsegment
%     \move (0 -15) \lvec (0 0) \onedot
% %     \lvec (-32 45)
%     \move (0 0) \lvec (-19 55)
%     \move (0 0)
%     \lvec (3 20) \onedot \lvec (-5 58)
%     \move (0 0) \lvec (20 15) \onedot
%     \lvec (12 30) \onedot %\lvec (9 55)
%     \move (20 15) 
%     \lvec (30 30) \onedot \rlvec (-4 21) \rmove (4 -21)
%     \rlvec (6 18)
%     \move (20 15) \rlvec (0 38)
%   \esegment
\end{texdraw}
\end{center}

\end{blanko}

One feature of operadic trees is that they admit colimit 
descriptions of basic operations: most importantly, grafting can be expressed 
as a pushout, and every tree is the colimit of its elementary subtrees, 
i.e.~trees without inner edges \cite{Kock:0807}.  This makes them well suited for constructions
(rather than just decomposition).  Another advantage is that symmetries of the
original nested structure (for example a nesting of Feynman graphs) are better
captured by operadic trees than by combinatorial trees, as the above picture
also illustrates: the combinatorial tree has a symmetry which does not
reflect a symmetry in the nesting, and fails to detect the inner symmetries
of the nesting.  (The symmetry issues play a crucial role in the treatment of
Green functions, where the operadic viewpoint seems important~\cite{Galvez-Kock-Tonks:FdB}.)

% \bigskip
% 
% To enhance the tree, one can decorate the nodes by substructures from the
% original picture, but then, depending on the theory, not all decorations make 
% sense: a node with a certain decoration may admit only certain children,
% and only a certain number.  In the end it becomes necessary to decorate also
% the edges, in order to express these compatibilities, and to keep track
% of the possible inputs to a given node, by drawing small vacant incoming edges
% of each node.  All this leads to the notion of operadic tree.

\begin{blanko}{The bialgebra of operadic trees.}\label{comultmy}
  A {\em cut} of an operadic tree is defined to be a subtree containing the root
  --- note that the arrows in the category of operadic trees 
  (Section~\ref{Sec:trees})
  are arity preserving (\ref{sub}), meaning that if
  a node is in the subtree, then so are all the incident edges of that node.
  
  If $c:R\subset T$ is a subtree containing the root, then each leaf $e$ of
  $R$ determines an {\em ideal subtree} of $T$ (\ref{sub}), namely consisting of $e$ (which
  becomes the new root) and all the descendant edges and nodes.  This is still
  true when $e$ is also a leaf of $T$: in this case, the ideal tree is the
  trivial tree consisting solely of $e$.  Figuratively, this means
  that cuts can go through the leaves but are not allowed to go {\em above}
  the leaves.  Note also that the root
  edge is a subtree;  the ideal tree of the root edge is of course the tree
  itself.  This is the analogue of the cut-below-the-root in the combinatorial
  case.  For a cut $c:R\subset T$, define $P_c$ to be the forest consisting
  of all the ideal trees generated by the leaves of $R$.
  
  Let $\mybialg$ be the polynomial $\ground$-algebra on the set of isomorphism 
  classes of operadic trees.  With comultiplication defined on the generators
  by
  \begin{eqnarray*}
  \Delta:  \mybialg & \longrightarrow & \ \ \mybialg \tensor_\ground \mybialg  \\
    T & \longmapsto & \sum_{c:R\subset T}\!\! P_c \tensor R ,
  \end{eqnarray*}
  as for combinatorial trees, $\mybialg$ becomes a graded bialgebra.  It is not
  connected: $\mybialg_0$ is spanned by powers of the trivial tree
  $\inlineDotlessTree$ (including the empty power, which is the algebra unit $1$).
  These are all group-like, so one
  could obtain a connected bialgebra by imposing the equation $1 =
  \inlineDotlessTree$.  Note that the comultiplication for operadic trees is
  linear in the right-hand factor, unlike the comultiplication for combinatorial
  trees, where $\Delta(T)$ always contains a factor $T \tensor 1$ corresponding
  to the empty-cut-below-the-root.
  
  As an example, here are the 11 possible cuts of the tree in the previous
  picture:
   \begin{center}
  \includegraphics[scale=0.8, clip = true , bb = 0 6 380 260]{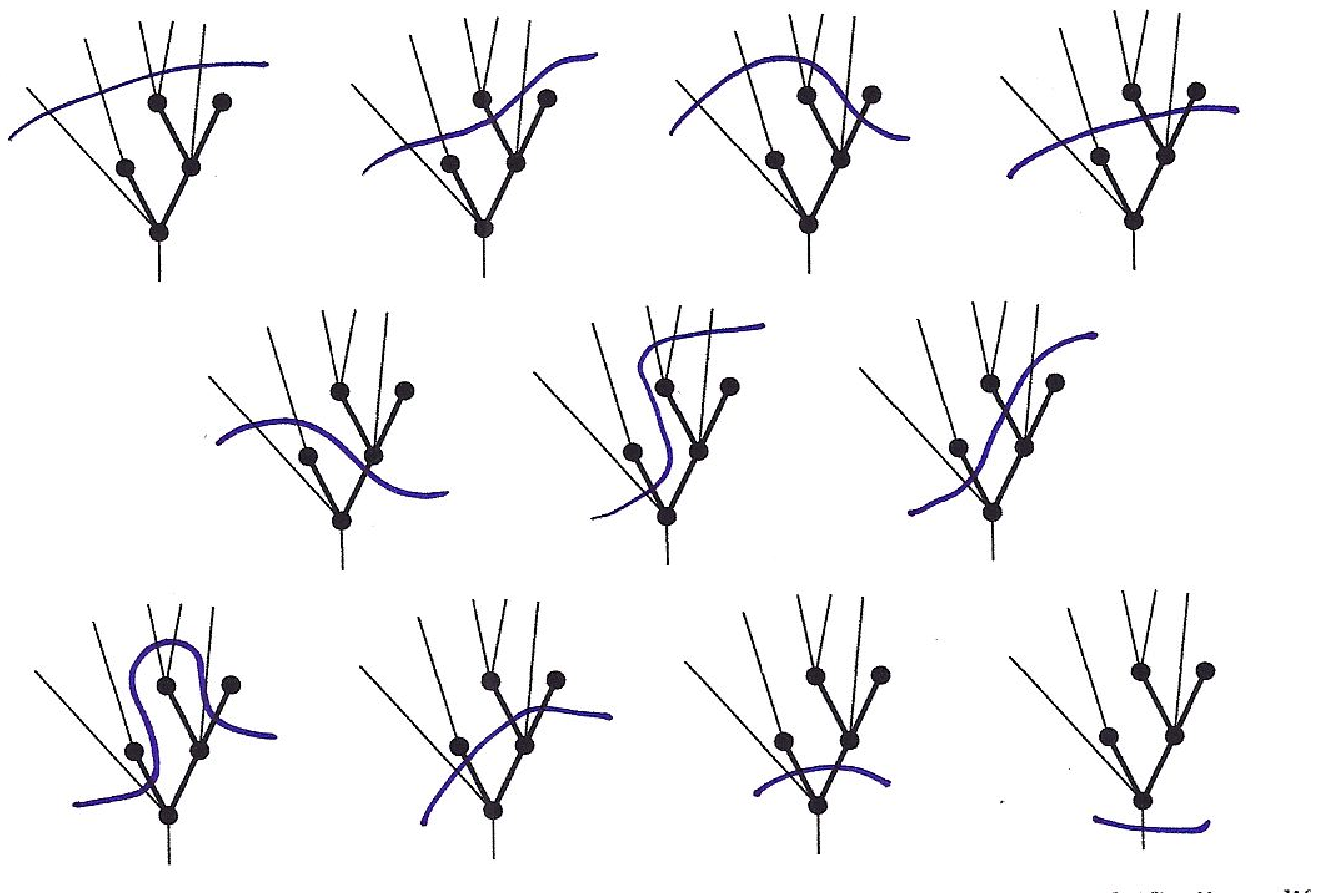}
  \end{center}

  Note that the first and the last term of the comultiplication of this tree $T$
  are
  $$
  \!\inlineDotlessTree\!\!\inlineDotlessTree\!\!
  \inlineDotlessTree\!\!\inlineDotlessTree\!\!\inlineDotlessTree\!
  \tensor T \ + \ T \tensor \!\inlineDotlessTree ,
  $$
  not $1\tensor T + T \tensor 1$. 

\end{blanko}
  
\begin{blanko}{From operadic trees to combinatorial trees.}
  The {\em core} of a non-trivial operadic tree $T$ is the combinatorial tree
  $T\upperdot$ obtained by pruning off all leaves as well as the root edge.
  Taking core is functorial in root-preserving inclusions.  For any non-trivial
  operadic tree $T$, there is a natural bijection between the set of non-trivial
  subtrees of $T$ and the set of combinatorial subtrees of $T\upperdot$.  This
  bijection sends $R \subset T$ to $R\upperdot \subset T\upperdot$.  By
  extending the assignment by defining the core of the trivial tree to be the
  empty (combinatorial) forest, it is clear that the operation is compatible
  with comultiplication, as illustrated in the previous picture,  where
  the core is highlighted with
  fat lines.
  
  In  conclusion we have:
\end{blanko}
  
\begin{blanko}{Proposition.}
  {\em
  Taking core defines a surjective homomorphism of graded bialgebras
  \begin{eqnarray*}
    \mybialg & \longrightarrow & \CK  \\
    T & \longmapsto & T\upperdot .
  \end{eqnarray*}}
\end{blanko}

\begin{blanko}{The bialgebra of $P$-trees.}
  The definition makes sense equally well for $P$-trees (\ref{Ptree}).  If there
  are more than one edge colour involved, there is a trivial tree for each edge
  colour, and the monomials in these trees (including the trivial monomial $1$)
  span $\mybialg_0$.  Clearly they are all group-like.
\end{blanko}

%%%%%%%%%%%%%%%%%%%%%%%%%%%%%%%%%%%%%%%%%%%%%%%%%%
\section{Polynomial functors}
%%%%%%%%%%%%%%%%%%%%%%%%%%%%%%%%%%%%%%%%%%%%%%%%%%

\label{Sec:poly}

In this section, we recall the basic notions of polynomial functors, 
referring to Gambino--Kock~\cite{Gambino-Kock:0906.4931} for all details.

\begin{blanko}{Slices, pullback, and adjoints.}
  We write $+$ and $\sum$ for the disjoint union of sets (i.e.~the categorical
  coproduct).  The equality sign denotes canonical isomorphism of sets.
  Let $B$ be a set.  Recall (e.g.~from~\cite{Awodey:Cat}, \S9.7) the {\em
  slice category} $\Set/B$: its objects are maps $X \to B$, and its arrows
  are commutative triangles
  \begin{diagram}[w=2ex,h=3.5ex,tight]
  X    && \rTo    && X'    \\
  &\rdTo    &      & \ldTo  &  \\
  &    & B .   & &
  \end{diagram}  
  For a set map $p:X \to B$, we denote the fibre over $b$ by $X_b
  := p^{-1}(b)$.  We can then write $X$ as the union of the fibres, $X =
  \sum_{b\in B} X_b$, and interpret $p$ as a $B$-indexed family of sets.
  In fact, there is a canonical equivalence of categories
  \begin{eqnarray*}
    \Set/B & \longrightarrow & \Set^B  \\
    {}[X \to B] & \longmapsto & [b\mapsto X_b] .
  \end{eqnarray*}
  Here $\Set^B$ denotes the category of functors $B \to \Set$ (considering
  the set $B$ as a discrete category).  For either interpretation of $B$-indexed
  family of sets, we shall
  use the notation $(X_b \mid b\in B)$.

  Given a map $f: A \to B$, we have the {\em lowershriek} functor
  \begin{eqnarray*}
   f\lowershriek:  \Set/A & \longrightarrow & \Set/B  \\
    {}[X\to A] & \longmapsto & [X\to A \to B].
  \end{eqnarray*}
  Using decomposition into fibres, the functor has the following description:
  $$
  (X_a\mid a\in A) \longmapsto (\sum_{a\in A_b} X_a \mid b\in B) ,
  $$
  i.e.~sum along the fibres.
  
  The functor $f\lowershriek$ has a right adjoint, denoted by {\em upperstar}, 
  given by pullback:
  \begin{eqnarray*}
    f\upperstar : \Set/B & \longrightarrow & \Set/A  \\
    {}[X\to B] & \longmapsto & [A\times_B X \to A] .
  \end{eqnarray*}
  In fibre notation,
  $$
  (X_b \mid b\in B) \longmapsto ( X_{f(a)} \mid a\in A ).
  $$
  
  Finally, also $f\upperstar $ in turn  has a right adjoint, denoted
  {\em lowerstar}, which is more involved to describe synthetically,
  but whose description in terms of fibres is ``multiply along the fibres'':
  \begin{eqnarray*}
    f\lowerstar :\Set/A  & \longrightarrow & \Set/B  \\
    (X_a\mid a\in A) & \longmapsto & (\prod_{a\in A_b} X_a \mid b\in B) .
  \end{eqnarray*}
  
  If the categories $\Set/A$ and $\Set/B$ are replaced by the equivalent 
  categories $\Set^A$ and $\Set^B$, the three functors $f\lowershriek$,
  $f\upperstar $, and $f\lowerstar $ can still be interpreted: $f\upperstar $
  is just precomposition with $f$.  Its adjoints are left and right Kan 
  extension respectively (cf.~\cite{MacLane:categories}, Ch.~X, for this 
  notion).  In this paper we shall not need the explicit form of these
  functors.  Although we actually formulate many results in terms of the
  functor categories $\Set^B$, when it comes to functors we prefer to work 
  with slices.
\end{blanko}

\begin{blanko}{Polynomial functors.}
  A diagram of sets
\begin{equation}\label{bridge}
  \begin{diagram}[w=3ex,h=3ex,tight]
  && E & \rTo^p  && B && \\
  & \ldTo^{s}&&&&&\rdTo^{t} & \\
I &&&&&&& J
\end{diagram}
\end{equation}
defines a {\em polynomial functor}
$$
\Set/I  \stackrel{s\upperstar}\rTo  \Set/E \stackrel{p\lowerstar}\rTo \Set/B 
\stackrel{t\lowershriek}\rTo  \Set/J .
$$
% Here upperstar, lowershriek and lowerstar denote pullback and its adjoints.
% In detail, 
%   for an arbitrary set map $f: A \to B$ we have the pullback functor $f\upperstar :
% \Set/B \to \Set/A$ and its adjoints:
% the left adjoint is $f\lowershriek:\Set/A\to\Set/B$ defined by composition, 
% i.e.~summing along the fibres, 
% $(X_a \mid a\in A) \mapsto ( \sum_{a\in A_b} X_a \mid b\in B)$, and
% the right adjoint is $f\lowerstar : \Set/A \to \Set/B$ whose description is
% $(X_a \mid a\in A) \mapsto ( \prod_{a\in A_b} X_a \mid b\in B)$, 
% i.e.~multiplying along the fibres. 
In terms of indexed families and fibres, the formula for this
polynomial functor is
\begin{eqnarray*}
  \Set/I & \longrightarrow & \Set/J  \\
  (X_i \mid i\in I) & \longmapsto & \big(\sum_{b\in B_j} \prod_{e\in E_b} 
  X_{s(e)} \mid j \in J \big) .
\end{eqnarray*}
It is hence an $J$-indexed family of sums of products of $I$-indexed families.
Particularly suggestive is the case $I=J=1$, so that we are talking about
a `single polynomial in one variable'.  In this case the formula boils down to
\begin{eqnarray*}
  \Set & \longrightarrow & \Set  \\
  X & \longmapsto & \sum_{b\in B} X^{E_b} .
\end{eqnarray*}

\end{blanko}

\begin{blanko}{Morphisms.}\label{nat}
  Morphisms of polynomial functors are just natural transformations.
  One can show (\cite{Gambino-Kock:0906.4931}) that a natural transformation 
  $P'\Rightarrow P$ between polynomial functors
  is uniquely
  represented by diagrams of the form
    \begin{equation}\label{equ:morphism}
  \begin{diagram}
  P': && I  &\lTo & E'& \rTo & B' &\rTo  &J \\
   && && \uTo && \dLig && \\
  &&\dLig &&  \bullet \SEpbk &  \rTo  & B' && \dLig \\
  && && \dTo && \dTo &&\\
   P : && I  &\lTo & E& \rTo & B &\rTo  &J . 
  \end{diagram}
  \end{equation}
  By $\Poly(I,J)$ we denote the category of polynomial functors $\Set/I \to 
  \Set/J$, and their natural transformations.
\end{blanko}

For the manipulation of polynomial functors in terms of the representing sets, 
the following two facts are basic.
\begin{blanko}{Beck--Chevalley.}
  Given a pullback square
  \begin{diagram}[w=4.5ex,h=4.5ex,tight]
  \ov A\SEpbk & \rTo^{\ov\fy}  & \ov B  \\
  \dTo<\alpha  &    & \dTo>\beta  \\
  A  & \rTo_{\fy}  & B
  \end{diagram}
  there are natural isomorphisms of functors
  $$
  \alpha\lowershriek \circ \ov\fy\upperstar  \isopil
  \fy\upperstar \circ \beta\lowershriek 
  \qquad \text{ and } \qquad 
    \beta\upperstar \circ \fy\lowerstar  \isopil \ov\fy{}\lowerstar \circ 
    \alpha\upperstar ,
    $$
   usually called the {\em Beck--Chevalley isomorphisms}.
\end{blanko}

\begin{blanko}{Distributivity.}
  Starting from maps $A \rTo^\fy B \rTo^\psi C$, we can construct the following
  diagram by applying $\psi\lowerstar $ to the map $\fy: A\to B$:
  \begin{diagram}[w=6ex,h=4.5ex,tight]
  \psi\upperstar\psi\lowerstar A\SEpbk & \rTo^{\ov\psi} & \psi\lowerstar A\\
  \dTo<{\epsilon} && \\
  A && \dTo>{\wtil{\fy}} \\
     \dTo<\fy   &  &   \\
    B & \rTo_\psi & C .
  \end{diagram}
  Here $\epsilon$ is the $A$-component of the  counit for the adjunction $\psi\upperstar 
  \isleftadjointto \psi\lowerstar $.  A diagram of this form is called a 
  {\em distributivity pentagon}; it can be characterised by a universal 
  property~\cite{Weber:1106.1983}.
  For such a diagrams the {\em distributive law} holds:
    $$
    \psi\lowerstar \circ \fy\lowershriek \simeq \wtil\fy\lowershriek 
    \circ \ov{\psi}{}\lowerstar \circ \epsilon\upperstar  .
    $$
  This is the categorical expression of the distributive
  law of elementary arithmetic, as it amounts to distributing a product 
  (lowerstar) over a sum (lowershriek);
  see~\cite{Gambino-Kock:0906.4931} for more discussion, and \cite{Weber:1106.1983}
  for a deeper treatment.
\end{blanko}

\begin{blanko}{Composition.}\label{comp}
  The composition of two polynomial functors is again poly\-nomial~\cite{Gambino-Kock:0906.4931}.
  This is a consequence of the Beck--Chevalley isomorphisms and distributivity.
  The important fact is that the bridge diagram representing the
  composite can be constructed explicitly in terms of a few operations on sets.
  Namely, given polynomial functors $P$ and $Q$ as in the bottom of the diagram
  (in which the labels $\Delta$, $\Pi$, and $\Sigma$ merely indicates which sort of
polynomial operation is performed along each map: $\Delta$ indicates pullback, 
$\Sigma$ indicates lowershriek, and $\Pi$ indicates lowerstar)
  \begin{diagram}[w=4.5ex,h=4.5ex,tight,labelstyle=\scriptscriptstyle]
  &&&&&&\cdot&&\rTo^\Pi&\cdot&\rTo^\Pi&\cdot\\
  &&&&&\ldTo^\Delta&\text{ \ \ \ \ \ pb}&&\ldTo_\Delta\\
&&&&\cdot&\rTo_\Pi&&\cdot&&\text{ \ \ \ \ \ distr}&&\dTo>\Sigma \\
&&&\ldTo^\Delta&\text{ \ \ \ \ \ pb}&&\ldTo_\Delta&&\rdTo_\Sigma&&&\\
&&\cdot&&\rTo_\Pi&\cdot&&\text{pb}&&\cdot&\rTo_\Pi&\cdot\\
&\ldTo^\Delta&&&&&\rdTo_\Sigma&&\ldTo_\Delta&&&&\rdTo^\Sigma\\
\cdot&&&Q&&&&\cdot&&&P&&&\cdot
\end{diagram}
the composite $P\circ Q$ is constructed as the top outline: start by taking the pullback at
the common middle set, then lowerstar the result to arrive at what will be the 
top right-hand corner of the composite, and pull back to complete the 
distributivity pentagon.  The diagram is completed by taking two more pullbacks 
as indicated.  The bridge diagram constructed is naturally isomorphic to the 
composite of the original functors by the Beck--Chevalley isomorphisms and 
distributivity.  The construction reflects closely how one composes two 
polynomial  functions in elementary algebra: the key point is of course 
distributing
the products involved in the outer polynomial over the sums of the inner.
\end{blanko}

\begin{blanko}{The $2$-category of polynomial functors.}
  Polynomial functors form a $2$-category\footnote{The notion of $2$-category 
  (see for example \cite{MacLane:categories}, Ch.~XII)
  is not essential to understand this work, but it is an efficient framework to 
  set up some of the involved notions correctly.} $\Poly$ in which the objects are the
  slices of $\Set$ (or equivalently, the categories $\Set^I$ for $I$ a set), the
  $1$-cells are the polynomial functors (i.e.~those isomorphic to one given by a
  diagram \eqref{bridge}), and $2$-cells are arbitrary natural transformations.
  $\Poly$ is a strict $2$-category.  By Theorem~2.17 of
  \cite{Gambino-Kock:0906.4931} it is biequivalent to a bicategory whose objects
  are sets $I$, whose $1$-cells are diagrams like \eqref{bridge} and whose
  $2$-cells are diagrams like \eqref{nat} (and where composition is described as
  in \ref{comp}).  The theorem allows us to blur the distinction between bridge
  diagrams and the polynomial functors they represent, allowing for the
  conceptual benefit of the strict $2$-category $\Poly$ and the computational
  benefit of the representing diagrams.  This is a characteristic aspect of the
  theory of polynomial functors.
\end{blanko}

\section{Elementary algebra of polynomial functors}
%%%%%%%%%%%%%%%%%%%%%%%%%%%%%%%%%%%%%%%%%%%%%%%%%%
%  \\[6pt]
% \normalsize{and na\"\i ve algebraic geometry over 
% the category of finite sets}

\label{Sec:alg}

We are here concerned with the aspect of polynomial functors as a categorification
of the elementary algebra of polynomial functions.  In this section we introduce a
set-up to deal with `rings' of polynomial functors and their associated `affine
spaces', in the sense of algebraic geometry~\cite{Hartshorne}.
The results in this section appear for the first time, but they are
not difficult.
One key issue is to impose the correct finiteness conditions.
Certainly, our `ground category' should be the category of finite sets: all
coefficients and exponents will now be required to be finite.  However, the
polynomial rings we are interested in have infinitely many variables, so we
cannot just take the theory of polynomial functors internally to the category of
finite sets.

\begin{blanko}{Finite sets.}
  The category of finite sets, which we denote
  $$
  \F := \FinSet,
  $$
  will be our coefficient `ring'.  More precisely, the monoidal operations of
  finite sums and products make it a distributive
  category~\cite{Carboni-Lack-Walters}.  Clearly the set of isomorphism classes
  of $\F$ is the set of natural numbers, and sum and products then yield the
  usual addition and multiplication of numbers.  In short,
  \begin{quote}
  {\em $\N$ is the Burnside semiring of $\F$},
  \end{quote}
  which is categorification in its purest form.
  
  The notation $\F$ stresses the algebraic aspect of $\FinSet$, and we 
  use this notation when we think of that category as the ground ring.
  The same 
  category also plays a geometric role, and as such we denote it
  $$
  \A :=\FinSet ,
  $$
  leading to a natural notation for what plays the role of affine space.
  Namely, if $I$ is a set, then the functor category
  $$
  \A^I = \FinSet^I
  $$
  is the domain for the finite polynomial functors in $I$-many variables.
  (Note that $\FinSet^I$ is not equivalent to $\FinSet/I$ when $I$ is an 
  infinite set: the latter is the category of maps $E \to I$ with $E$ finite.
  The former is equivalent to the category of maps $E \to I$ with finite fibres.)
\end{blanko}

\begin{blanko}{Finite polynomials.}
  A polynomial 
  $$
  I  \leftarrow E  \to B \to 1
  $$
  is called {\em finite} when $B$ and $E$ are finite sets. 
  These are the polynomial functors on $\A^I$ with values in finite sets.
  We denote the category of these polynomial functors
  $$
  \F[I] := \kat{FinPoly}(I,1) .
  $$
  Clearly $\F[I]$ is a distributive category and
  \begin{quote}
    {\em the Burnside semiring of $\F[I]$ is $\N[I]$, the polynomial semiring 
    in $I$-many variables,}
  \end{quote}
  as the notation also suggests. Moreover, 
  \begin{quote}
    {\em $\F[I]$ is the free distributive category on the set $I$,}
  \end{quote}
  just as $\N[I]$ is the free commutative semiring on $I$.  Indeed, the category
  $\F[I]$ is generated freely under finite sums by the monomials, i.e.~those $I
  \leftarrow E \to B \to 1$ for which $B=1$.  From the characterisation of
  natural transformations in \ref{nat}, we see that the category of monomial
  functors in $I$-many variables has arrows given by commutative triangles
  \begin{diagram}[w=4ex,h=3.5ex,tight]
   &   & E'  \\
  I  &  \ldTo(2,1) & \uTo  \\
    &\luTo(2,1) & E ,
  \end{diagram}
  and so is equivalent to $(\FinSet/I)\op$.  But it is well known that
  $(\FinSet/I)\op$ is the finite-product completion of $I$.
\end{blanko}

\begin{blanko}{`Polynomial rings' and `affine space'.}
  The distributive category $\F[I]$ is the category of finite
  polynomial functors on $\A^I$.  Conversely, the category $\A^I$ can be
  reconstructed from $\F[I]$: define a
  {\em character} on a distributive category $\DD$ to be a functor $\DD \to \F$
  preserving finite sums and finite products. The form a category in which the 
  morphisms are the finite-sum-and products-compatible natural transformations.
  \begin{quote}
    {\em $\A^I$ is equivalent to the category of characters of $\F[I]$.}
  \end{quote}
  Indeed, the category of sum-and-product
  preserving functors $\F[I]\to\F$ is equivalent to the product-preserving
  functors $(\FinSet/I)\op \to \F$, and since the domain is the 
  product-completion of $I$, we are finally left with the category of
  functors $I \to \F$, which is what we denote $\A^I$.
\end{blanko}

\begin{blanko}{Appropriate maps: locally finite polynomial functors.}
  \label{duality}
  For a given set $I$ we now have the category $\F[I]$ playing the role of a
  polynomial ring, and the category $\A^I$ playing the role of affine space.  We
  proceed to assemble these objects into $2$-categories $\kat{PolyAlg}$ and
  $\kat{Aff}$, respectively, which will be the ambient setting for defining
  the algebraic structures promised in the introduction.
%   .  The $2$-category $\kat{PolyAlg}$ will be the
%   setting for defining tensor products like $\F[I]\tensor \F[I]$, needed in
%   order to say what a comultiplication is.
  
  The $2$-category $\kat{PolyAlg}$ is defined as follows.  Its objects are the
  categories of the form $\F[I]$ (i.e.~free distributive categories on a set).
  The $1$-cells are the finite-sum-and-product-preserving functors, and the
  $2$-cells are the compatible natural transformations (i.e.~those whose 
  component on a finite sum is the sum of the components, and similarly with 
  finite products).  The universal property of $\F[J]$ implies that
  the hom cats in $\kat{PolyAlg}$ are
  $$
  \kat{PolyAlg}( \F[J],\F[I]) \simeq \F[I]^J .
  $$
  
%   Given two categories of finite polynomial functors, $\F[J]$ and $\F[I]$, the
%   appropriate category of functors from $\F[J]$ to $\F[I]$ consists of
%   finite-sum-and-product-preserving functors and compatible natural
%   transformations.
  
  On the other hand, the $2$-category $\kat{Aff}$ is defined as having
  the categories $\A^I$ as objects, and as hom cats the categories of
  locally finite polynomial functors and their natural transformations,
  denoted $\kat{LocFinPoly}(I,J)$:  
  A polynomial functor
  $$
  I \stackrel s \leftarrow E \stackrel p \to B \stackrel t \to J
  $$
  is called {\em locally finite} when the maps $p$ and $t$ have finite fibres.
  This means that only finite sums and finite products are involved.
  It is not difficult to see that we have an equivalence of categories
  $$
  \kat{Aff}(\A^I,\A^J) := \kat{LocFinPoly}(I,J) \simeq \F[I]^J .
  $$
  One may observe that the locally finite polynomial functors can also be
  charaterised as those
  $F: \A^I \to \A^J$ such that for any $P$ belonging to $\F[J]$, the composite
  $P\circ F$ belongs to $\F[I]$, in analogy with the definition of regular 
  maps in algebraic geometry.
%   The universal property of $\F[J]$ readily implies that in both cases the
%   functor category is $\F[I]^J$.  
%   
%   The category $\F[I]^J$ has the following expected 
%   explicit description.

%   From the explicit construction of composition in \ref{comp} it follows
%   readily that the composite of two locally finite polynomial functors is again
%   locally finite.  Altogether we have a $2$-category $\kat{Aff}$ whose objects
%   are the categories of the form $\A^I$, and whose hom categories are
%   the $\kat{LocFinPoly}(I,J)$.  On the other hand, we have the 
%   $2$-category $\FinPoly$ of categories of finite polynomial functors,
%   finite-sum-and-product-preserving functors and compatible natural
%   transformations. 
  
  Almost tautologically we have an equivalence of
  $2$-categories
  $$
  \kat{Aff} \simeq \kat{PolyAlg}\op ,
  $$
  justifying the symbols and the terminology.
\end{blanko}

\begin{blanko}{The tensor product and comonoidal structure.}
  The $2$-category $\Aff$ has finite products: it is given by the equivalence
  of categories
  $$
  \A^{I_1} \times \A^{I_2}\simeq \A^{I_1+I_2} .
  $$
  Accordingly, the category $\kat{PolyAlg}$ has categorical sums, which we denote by
  $\tensor_\F$.  
%   This is a special case of the Gates tensor product of
%   distributive categories.
  All we need to know about it is the equivalence
$$
\F[I_1] \tensor_{\F} \F[I_2] \simeq  \F[I_1+I_2] ,
$$
which is just dual to the previous display.  The neutral element for this tensor
product is the category $\F=\F[\emptyset]$ itself, which is an initial object in
$\kat{PolyAlg}$: the unique sum-and-product-preserving functor $\F\to\F[I]$
sends a finite set $S$ to the polynomial functor
$$
I \leftarrow 0 \to S \to 1 .
$$
% The tensor product on $\kat{PolyAlg}$ is the restriction of the so-called Gates 
% tensor product of distributive categories \cite{Gates}, referred to by Lawvere and 
% Menni~\cite{Lawvere-Menni}.

The objects in $(\kat{PolyAlg}, \tensor, \F)$ are algebras.  We would like to 
define bialgebras as comonoids in $(\kat{PolyAlg}, \tensor, \F)$, or rather 
pseudo-comonoids (i.e.~coassociative only up to given coherent isomorphisms).
A bialgebra is then a category $\F[I]$ equipped with a comonoidal structure,
i.e.~functors $\F \leftarrow \F[I] \to \F[I]\tensor \F[I]$.  Every $\F[I]$
has two canonical such comonoidal structures, as we shall see in the next 
paragraph.
The main result of this paper establishes another comonoidal structure
for the special case $I=\tr$.
However, it is much more convenient to describe these structures on the other
side of the duality, namely in $(\kat{Aff}, \times, 1)$, where the background
structure is just the product, and we are talking about the familiar
notion of
monoidal structure instead of comonoidal structure.
\end{blanko}

% 
% \bigskip
% 
% 
% 
% The distributive category $\F[I]$ has a natural sum-and-product-preserving
% functor from $\F$ assigning to a set $S$ the constant polynomial with 
% value $S$; this is represented by
% $$
% I \leftarrow 0 \to S \to 1 .
% $$
% This is to say that $\F[I]$ is an $\F$-algebra.
% It is the free $\F$-algebra on the set $I$ (perhaps in a suitable 
% $2$-categorical sense).  The categorical sum
% in the category of distributive categories
% is the tensor product $\tensor_{\F}$, which we need to work out.
% Since the functor $\Set \to \DistrCat$ is a left adjoint (in some 
% $2$-categorical sense?), it preserve sums, in particular we find
% $$
% \F[I_1] \tensor_{\F} \F[I_2] \simeq  \F[I_1+I_2] .
% $$

\begin{blanko}{Canonical monoidal structures on $\A^I$.}
  Fix a set $I$.  The category $\A^I$ has two canonical monoidal structures:
  pointwise sum and pointwise product.  These are given by locally finite
  polynomial functors.
  
  The {\em sum} of $x:I\to\A$ and $y:I\to\A$ is 
  the function $x+y: I \to \A$ which sends $i$ to $x(i)+y(i)$.
  Under the equivalence
  $$
  \A^I \times \A^I \simeq \A^{I+I}
  $$
  the sum map is polynomial: it is represented by
  $$
  I+I \lTo I+I \rTo I+I \rTo I .
  $$
  The corresponding comonoidal structure on $\F[I]$ is given by
  \begin{eqnarray*}
    \Delta:\F[I]& \longrightarrow & \F[I]\tensor_\F \F[I]   \\
    x_i  & \longmapsto & x_i\tensor 1 + 1 \tensor x_i.
  \end{eqnarray*}
  Here $x_i \in \F[I]$ denotes the polynomial functor
  $$
  I \stackrel{\name{i}}\lTo 1 \stackrel = \rTo 1 \stackrel = \rTo 1 .
  $$
  where $\name i : 1 \to I$ picks out the element $i$.
  
  The neutral object for the sum is of course the constant empty set $0$.
  Under the equivalence $1 \simeq \A^0$,
  the corresponding map is represented by
  $$
  0 \leftarrow 0 \to 0 \to I.
  $$
  The corresponding counit for the comultiplication is given by $\epsilon(x_i)=0$.
  
  The {\em product} of $x:I\to\A$ and $y:I\to\A$  is the function
  $xy: I \to \A$ which sends $i$ to $x(i)y(i)$.  It is represented by
  $$
  I+I \lTo I+I \rTo I \rTo I .
  $$
  The unit is the constant function $1$.  The corresponding comonoidal structure
  on $\F[I]$ is given by $\Delta(x_i) = x_i\tensor x_i$ and $\epsilon(x_i)=1$.
  This categorifies the construction of
  the monoid-algebra on the
  free abelian monoid on $I$.
\end{blanko}

Before coming to the promised new monoidal structure on $\A^\tr$ in 
Section~\ref{Sec:monoidal}, we need some background on trees.

\section{$P$-trees and free monads}
%%%%%%%%%%%%%%%%%%%%%%%%%%%%%%%%%%%%%%%%%%%%%%%%%%

\label{Sec:trees}

We briefly summarise the polynomial formalism of trees from \cite{Kock:0807}.

\begin{blanko}{Trees.}\label{polytree-def}
  It was observed in \cite{Kock:0807} that operadic trees can be
  conveniently encoded by diagrams of the same shape as polynomial functors.
  We take this as the definition of tree: An
  {\em operadic tree} is a diagram of finite sets
    \begin{diagram}
    A & \lTo^s & M & \rTo^p & N & \rTo^t & A
\end{diagram}
satisfying the following three conditions:
  
  (1) $t$ is injective
  
  (2) $s$ is injective with singleton complement (called the {\em 
  root} and denoted $1$).
  
  \noindent With $A=1+M$, 
  define the walk-to-the-root function
  $\sigma: A \to A$ by $1\mapsto 1$ and $e\mapsto t(p(e))$ for
  $e\in M$. 
  
  (3)  $\forall x\in A : \exists k\in \N : \sigma^{k}(x)=1$.
  
  The elements of $A$ are called {\em edges}.  The elements of $N$
  are called {\em nodes}.  For $b\in N$, the edge $t(b)$ is called
  the {\em output edge} of the node.  That $t$ is injective is just to
  say that each edge is the output edge of at most one node.  For
  $b\in N$, the elements of the fibre $M_b\df p^{-1}(b)$ are
  called {\em input edges} of $b$.  Hence the whole set
  $M=\sum_{b\in N} M_b$ can be thought of as the set of
  nodes-with-a-marked-input-edge, i.e.~pairs $(b,e)$ where $b$ is a
  node and $e$ is an input edge of $b$.  The map $s$ returns the
  marked edge.  Condition (2) says that every edge is the input edge
  of a unique node, except the root edge.
  Condition (3) says that if you walk towards the root, in a finite 
  number of steps you arrive there.

  The edges not in the image of $t$ are called {\em leaves}.
  From now on we just say {\em tree} for `operadic tree'.
  
  The tree
  $$
  1 \leftarrow 0 \to 0 \to 1
  $$
  is the {\em trivial tree}, which we denote by \inlineDotlessTree.
\end{blanko}

\begin{blanko}{Cartesian morphisms (cf.~\cite{Gambino-Kock:0906.4931}).}
  \label{cart}
  A {\em cartesian morphism} of polynomial endofunctors is by definition a diagram
  \begin{equation}
    \label{equ:cartmorphism}
  \begin{diagram}[w=4.5ex,h=3.5ex,tight]
  I' & \lTo & E' \SEpbk &\rTo & B'& \rTo & I'  \\
  \dTo<\alpha && \dTo && \dTo && \dTo>\alpha \\
  I  &\lTo & E& \rTo & B &\rTo   &I .
  \end{diagram}
  \end{equation}
  If $\alpha$ is an identity map, then this corresponds
  precisely to cartesian natural transformations between polynomial endofunctors
  (i.e.~whose naturality squares are cartesian).  In the general case, the
  cartesian morphisms like \eqref{equ:cartmorphism} are outside the scope of the
  $2$-category $\Poly$, but they can be reduced to it by some base
  change~\cite{Gambino-Kock:0906.4931}.
\end{blanko}

\begin{blanko}{Morphisms of trees (cf.~\cite{Kock:0807}, \S 1.1).}\label{sub}
  A {\em morphism of trees} is by definition a cartesian morphism of the 
  associated polynomial
  endofunctors.  Hence edges go to edges, and the cartesian
  condition means that a node of arity $n$ is mapped to a node of arity $n$.
  The morphisms of trees are also called {\em tree embeddings} since in
  fact each of
  the components of such a map is automatically injective: indeed, the 
  tree axioms above imply that a morphism necessarily commutes with the 
  walk-to-the-root function $\sigma$, and together with arity perservation this
  forces the maps to be injective.  Hence the category of
  trees and tree embeddings, denoted $\TEmb$, is mostly concerned with subtrees,
  but note that it also contains automorphisms of trees.
  
  A tree embedding is {\em root-preserving} when it sends the root to the root.
  In formal terms, these are diagrams \eqref{equ:cartmorphism}
  such that also the left-hand square is cartesian, \cite{Kock:0807}, 1.1.13.
  
  An {\em ideal embedding} (or an {\em ideal subtree}) is a subtree for which every
  descendent edge and node is also in the tree, \cite{Kock:0807}, 1.1.9.  There is one ideal subtree
  generated by each edge in the tree.  The ideal embeddings are characterised
  as having also the right-hand square of \eqref{equ:cartmorphism} cartesian.
  
  Ideal embeddings and root-preserving embeddings admit pushouts along each 
  other in the category $\TEmb$, \cite{Kock:0807}, 1.1.20.  The most interesting case is pushout over a 
  trivial tree: this is then the root of one tree and a leaf of another tree,
  and the pushout is the grafting onto that leaf.
\end{blanko}

\begin{blanko}{Decorated trees: $P$-trees.}\label{Ptree}
  Let $P$ denote a finitary polynomial endofunctor on $\Set/I$, represented by
  $$
  I \leftarrow E \to B \to I .
  $$
  By definition (\cite{Kock:0807}, \S 1.2), a {\em $P$-tree} is a diagram
  \begin{diagram}[w=4.5ex,h=3.5ex,tight] 
  A&\lTo&M \SEpbk&\rTo&N &\rTo& A\\
  \dTo && \dTo && \dTo && \dTo \\
  I&\lTo&E&\rTo&B &\rTo& I   .
  \end{diagram}
  where the top row is a tree.
  Unfolding the definition, we see that a $P$-tree is a
  tree whose edges are decorated in $I$, whose nodes are decorated
  in $B$, and with the additional structure of a bijection for each
  node $n \in N$ (with decoration $b \in B$) between the set of
  input edges of $n$ and the fibre $E_b$, subject to the
  compatibility condition that such an edge $x\in M_n$ corresponding to $e\in 
  E_b$ has
  decoration $s(e)$, and the output edge of $n$ has decoration $t(b)$.
  Note that the $I$-decoration of the edges is completely specified
  by the node decoration together with the compatibility requirement,
  except for the case of a trivial tree.
\end{blanko}

\begin{blanko}{The free monad on a polynomial endofunctor.}
  (For the notion of monad, see \cite{MacLane:categories}, Ch.VI, or 
  \cite{Awodey:Cat}, Ch.~10.)
  By a {\em polynomial monad} we mean a monad in the $2$-category of polynomial
  functors and cartesian natural transformations,
  cf.~\cite{Gambino-Kock:0906.4931}.  The forgetful functor from polynomial
  monads to polynomial endofunctors has a left adjoint, the free-monad functor.
  It has the following pleasing description in terms of trees, cf.~\cite{Kock:0807}, 
  Prop.~1.2.8:
  the free monad on a polynomial endofunctor $P$ is the polynomial monad
  represented by
  $$
  I \leftarrow \tr'_P \stackrel q \to \tr_P \to I
  $$
  where $\tr_P$ is the set of (isomorphism classes of) $P$-trees, and $\tr'_P$
  is the set of (isomorphism classes of) $P$-trees with a marked leaf.
  The left-most map returns the decoration of the marked leaf; the right-most 
  map
  returns the decoration of the root, and the middle map $q$ just forgets the 
  mark.
  The 
  monad
  structure is given by grafting of trees, an operation that allows an 
  easy and completely formal description in the category of trees in terms of 
  pushouts, cf.~\cite{Kock:0807}, Prop.~1.1.19.  Examples of the notion of $P$-trees, for
  suitable choices of $P$, are planar and binary trees, or more specialised 
  examples like trees with 
  nodes decorated by primitive 1PI graphs of
  a quantum field theory~\cite{Kock:graphs-and-trees},
  or the opetopes of higher category theory~\cite{zoom}.
  
  The polynomial endofunctor $P$ will be held fixed throughout, and from 
  now on we put $\tr:= \tr_P$.  We let $F$ denote the free monad on
  $P$.
  
  We shall need an explicit description of the monad structure of $F$.
  As in \ref{comp},
the composite $F\circ F$, is built by the following diagram, whose constituents 
we now make explicit:
\begin{diagram}[w=4.5ex,h=4.5ex,tight]
  &&&&&&\sect{\tr}{}'&&\rTo&\ptsect \tr&\rTo^p&\sect \tr&& = q\lowerstar (\tr\times_I \tr')\\
  &&&&&\ldTo&\text{ \ \ \ \ \ \footnotesize{pb}}&&\ldTo_\epsilon\\
&&&&\cdot&\rTo&&\tr\times_I \tr'&&\text{ \ \ \ \ \ \footnotesize{distr}}&&\dTo>r \\
&&&\ldTo&\text{ \ \ \ \ \ \footnotesize{pb}}&&\ldTo&&\rdTo&&&\\
&&\tr'&&\rTo_q&\tr&&\text{\footnotesize{pb}}&&\tr'&\rTo_q&\tr\\
&\ldTo&&&&&\rdTo&&\ldTo&&&&\rdTo\\
I&&&F&&&&I&&&F&&&I
\end{diagram}
The construction starts from the bottom by forming the central pullback square:
$\tr\times_I \tr'$ is the set of pairs $(S,R)$ such that the root of $S$
is of the same type as the marked leaf of $R$, in other words
a simple grafting of one tree onto another.  This can also 
be
interpreted as the set of trees with one marked edge: the two projections are
then: return the ideal subtree generated by that edge (that's $S$), and return the 
root-preserving tree obtained by pruning that ideal subtree (that's $R$).

Next we compute $\sect{\tr} := q\lowerstar (\tr\times_I \tr')$. (Note that this 
set can also be described as $F(\tr)$.)
By definition
of the lowerstar operation, the result is a set over $\tr$ whose fibre over a
fixed tree $R\in\tr$ is the set of maps from the set of leaves of $R$ to the set of all
trees with matching root.  This data is equivalent to the grafting of all those
leaf-indexed trees onto the leaves of $R$.  Alternatively this data amounts to
the inclusion of $R$ into the big tree resulting from the graftings.  In other
words, $\sect{\tr}$ is the set of trees with a cut, and the map $r: \sect{\tr} 
\to \tr$ returns the tree found below the cut.

The pullback along $q$ is the set $\ptsect{\tr}$ of pairs consisting of a marked tree and graftings
onto all leaves, which amounts to trees with a pointed cut; the map $p: 
\ptsect{\tr} \to \sect{\tr}$ just forgets the mark in the cut.
The evaluation map  $\epsilon:\ptsect{\tr}\to \tr\times_I \tr'$ is 
described like this: if we think
of the elements of $\ptsect{\tr}$ as given by a tree together with a marked leaf and a map
from each leaf to $\tr$ (with matching root), then the map $\epsilon$ simply 
applies that map to the marked leaf hence obtaining a single tree with correct 
root, and hence an element in $\tr\times_I \tr'$. 

(We don't really need the remaining left-hand part of the diagram, but here it 
is, for completeness: the upper left-hand corner is the set of trees with a pointed cut
and a marked leaf in the subtree corresponding to the point.  (Since that marked
leaf together with the cut automatically provides a point on the cut
it is enough to say: tree with a cut and a marked leaf.)  The set just
below it is the set of trees with a marked edge and a marked descendant leaf.)

Finally, we describe the structure maps of the free monad $F$:
the multiplication
$F\circ F \Rightarrow F$
is the cartesian natural
transformation
  \begin{diagram}[w=4.5ex,h=3.5ex,tight] 
  F\circ F:&& I&\lTo& \sect{\tr}{}'\SEpbk&\rTo&\sect\tr &\rTo& I\\
  \Downarrow \ \ &&\dLig && \dTo && \dTo>m && \dLig \\
 F: && I&\lTo&\tr'&\rTo&\tr &\rTo& I   .
  \end{diagram}
where the maps in the middle simply forget the cut.  
The unit for the monad, $\Id \Rightarrow F$ is 
given essentially by the map $e: I \to \tr$ assigning to each `colour' $i\in I$
the trivial tree with edge of colour $i$.

\bigskip

We summarise the maps $r$, $p$, $f$, $m$, for later use:
\begin{diagram}[w=4.5ex,h=4.5ex,tight]
  &&&&\ptsect{\tr}&\rTo^p&\sect{\tr}&=F(\tr)\\
  &\begin{picture}(0,0)(-20,16)
\qbezier(45,45)(-30,30)(-45,-45)
\put(-24,20){\makebox(0,0)[t]{$\scriptstyle f$}}
\put(-45,-45){\vector(-1,-4){0}}
\end{picture}&&\ldTo\\
&&\tr\times_I \tr'&&&&\dTo>r \\
&\ldTo&&\rdTo&&&&\begin{picture}(0,0)(6,16)
\qbezier(-12,87)(30,0)(-15,-95)
\put(16,10){\makebox(0,0)[t]{$\scriptstyle m$}}
\put(-14.5,-95){\vector(-1,-4){0}}
\end{picture}\\
\tr&&&&\tr'&\rTo&\tr\\
&\rdTo&&\ldTo&&&&\rdTo\\
&&I&&&F&&&I\\
&&&&&&&\ruTo\\
&&&&&&\tr
\end{diagram}
% So it is all about taking $F$ on the set of trees!  The maps appearing here
% are all for the left-hand part of the comultiplication diagram below.  The right-hand part
% of the comultiplication diagram will be the
% monad multiplication law $m: \sect{\tr}=F(\tr)=FF(1) \to F(1) = \tr$.

% [Question: is the free monad a coalgebra?  I.e.~can be construct a levelwise
% comultiplication by assigning to $F(X)$ ($X$-decorated trees) the sum all ways
% of putting a cut on all these trees?]
\end{blanko}

%%%%%%%%%%%%%%%%%%%%%%%%%%%%%%%%%%%%%%%%%%%%%%%%%%
\section{New monoidal structure on $\A^\tr$}
%%%%%%%%%%%%%%%%%%%%%%%%%%%%%%%%%%%%%%%%%%%%%%%%%%

\label{Sec:monoidal}

In this section the symbol $\tr$ stands for the set of
isomorphism classes of $P$-trees, for any fixed finitary polynomial endofunctor $P$
(and similarly for the decorated symbols, $\sect\tr$, etc.).
The construction of the monoidal structure works also for abstract trees,
but then the various sets and maps no longer come from the free-monad 
construction, and need to be defined in a more ad hoc manner.
We discuss this in the next section.

\begin{blanko}{Theorem.}
  {\em
  The polynomial functor $M:\A^\tr \times\A^\tr \to \A^\tr$ defined by
\begin{diagram}
  && \ptsect{\tr}+\sect{\tr} & \rTo^{\langle p,\sect{\tr}\rangle}  && \sect{\tr} && \\
  & \ldTo^{f+r}&&&&&\rdTo^{m} & \\
\tr+\tr &&&\ \ \ M&&&& \tr
\end{diagram}
is a monoidal structure on $\A^\tr$.  Its unit is  the functor $U:\A^0 \to \A^\tr$
given by
\begin{diagram}
  && 0 & \rTo  & I && \\
  & \ldTo&&&&\rdTo^{e} & \\
0 &&&U&&& \tr .
\end{diagram}
}\end{blanko}

\begin{blanko}{Description of the maps.}
An interesting feature of this monoidal structure is that all the involved maps
are recognised as those occurring in the polynomial description of the 
free-monad construction.
We review the definitions of the sets and maps:

\begin{itemize}
  \item $\tr$ is the set of iso-classes of $P$-trees.  
If $F$ denotes the free monad on $P$, then $\tr=F(1)$.

  \item The set $\sect{\tr}$ is the set of iso-classes of $P$-trees with a cut.
It appears as $F(\tr) = FF(1)$.  More formally, the set of iso-classes
of trees with a subtree 
containing the root.

  \item The map $m : \sect{\tr} \to \tr$ is the
multiplication of the monad, $m: FF(1) \to F(1)$, i.e.~forget the cut.

  \item 
The set $\ptsect{\tr}$ is the set of iso-classes of
$P$-trees with a pointed cut,
and $p$ is the map that forgets the point.
More formally,
$$
\ptsect{\tr} = \tr' \times_\tr \sect{\tr} = q\upperstar \sect{\tr} ,
$$
where $\tr'$ is the set of iso-classes of $P$-trees with a marked leaf (it appears as the 
total space of $q: \tr'\to \tr$ 
in the diagram representing $F$), and $p: \ptsect{\tr}\to \sect{\tr}$ is the projection. 

  \item The map
$r: \sect{\tr}\to \tr$ returns the $P$-tree below the cut. 
% More formally it is simply the structure map $F(\tr) \to \tr$.

\item $f: \ptsect{\tr}\to \tr$ returns the ideal subtree
generated by the edge marked by the pointed cut.  This is the most complicated
map to explain formally, but it occurs already in the diagram for $F\circ F$, and
is composed of the $\tr$-component of the counit of the $q\upperstar
\isleftadjointto q\lowerstar $ adjunction (that's the evaluation map 
$\epsilon:q\upperstar \sect{\tr}\to
\tr\times_I \tr'$) followed by the projection to $\tr$.

\item Finally, $e: I \to \tr$ assigns to each `colour' $i\in I$
the trivial tree with edge of colour $i$.
\end{itemize}

Note that here and throughout, the name of a set is also used to indicate its
identity map.  So for example ${\langle p,\sect{\tr}\rangle}$ denotes the map
whose first component is $p$ and whose second component is the identity map on
$\sect\tr$.

\bigskip

Note that the polynomial is a `sum-wise tensor product' of a non-linear functor 
(in the first set of variables):
\begin{diagram}
  && \ptsect{\tr} & \rTo^p  && \sect{\tr} && \\
  & \ldTo^{f}&&&&&\rdTo^{m} & \\
\tr &&&&&&& \tr
\end{diagram}
(which concerns taking all cuts and retaining the forest of cut-off branches), 
and a linear functor (in the second set
of variables):
\begin{diagram}
  && \sect{\tr} & \rTo^{=}  && \sect{\tr} && \\
  & \ldTo^{r}&&&&&\rdTo^{m} & \\
\tr &&&&&&& \tr
\end{diagram}
(which is about taking all cuts and retaining the bottom tree, which of course 
is a single tree, hence the linearity).
\end{blanko}

\begin{blanko}{Proof of associativity.}
  In the following diagram, which commutes strictly,
  the left-hand vertical polynomial is $\tr+M$, and the bottom is $M$.
  Similarly the top polynomial is $M + \tr$ and the right-hand vertical
  polynomial is $M$.  We show that both composites are naturally isomorphic to
  the polynomial indicated by the diagonal (to be detailed in the proof).  In
  both cases this amounts to performing the constructions as in \ref{comp}.
  This natural isomorphism provides the associator, which is part of the
  monoidal structure.  (The pentagon equation for the associator is not
  established explicitly.  We invoke instead a general coherence principle: since the
  associator is constructed by canonical isomorphisms (Beck--Chevalley and 
  distributivity), it is coherent.)
  
\begin{diagram}[w=40pt,h=30pt,tight,scriptlabels]
\tr\!+\!\tr\!+\!\tr&&\lTo^\Delta&&(\ptsect{\tr}\!+\!\sect{\tr})\!+\!\tr&\rTo^\Pi &\sect{\tr}+\tr&\rTo^\Sigma&\tr+\tr\\
&\luDashto(4,4)^\Delta&&&\uTo<\Delta&\text{\tiny pb}&\uTo>\Delta&\text{\tiny pb}&\uTo>\Delta\\
\uTo<\Delta&&&&\cdot&\rTo^\Pi&\cdot&\rTo^\Sigma&\ptsect{\tr}+\sect{\tr}\\
&&&&\uTo<\Delta&\text{\tiny pb}&\uTo>\Delta&&\\
\tr\!+\!(\ptsect{\tr}\!+\!\sect{\tr})&\lTo^\Delta&\cdot&\lTo^\Delta&\cdot&\rTo^\Pi&\cdot&\text{\tiny distr}&\dTo>\Pi\\
\dTo<\Pi&\text{\tiny pb}&\dTo<\Pi&\text{\tiny pb}&\dTo<\Pi&\rdDashto^\Pi&\dTo>\Pi&&\\
\tr+\sect{\tr}&\lTo_\Delta&\cdot&\lTo_\Delta&\cdot&\rTo_\Pi&\sectsect{\tr}&\rTo^\Sigma&\sect{\tr}\\
\dTo<\Sigma&\text{\tiny pb}&\dTo<\Sigma&&\text{\tiny distr}&&\dTo<\Sigma&\rdDashto^\Sigma&\dTo>\Sigma\\
\tr+\tr&\lTo_\Delta&\ptsect{\tr}+\sect{\tr}&&\rTo_\Pi&&\sect{\tr}&\rTo_\Sigma&\tr .
\end{diagram}

We start with the lower left-hand composite.
The first step is to take a pullback in the lower left-hand corner:
\begin{diagram}[w=7ex,h=4.5ex,tight]
\tr+\bluesect{\tr}&\lTo^{\red{f+r}}&\redptsect{\tr}+\redsectbluesect{\tr}  \\
\dTo<{\tr+\blue{m}}&\text{\tiny pb}&\dTo>{\redptsect{\tr}+\blue{\tilde m}}  \\
\tr+\tr&\lTo_{\red{f+r}}&\redptsect{\tr}+\redsect{\tr}.
\end{diagram}
Here $\redsectbluesect\tr$ is the set of trees with two non-crossing cuts, or 
more  formally, a sequence of two root-preserving inclusions $R\subset S \subset T$.
The upper cut is the one coming from the lower right-hand corner of the diagram,
and the lower cut comes from the upper left-hand corner of the diagram.
Hence the map denoted $\blue{\tilde m}$ forgets the lower cut.%
\footnote{Typographical note: if colour output is available, as an extra visual aid
the upper cut and the maps related to it
are printed in red, while the lower cut and its maps in blue.  The wording is,
however, intended to be sufficient for also a black-and-white printing to
make sense.}

The right-hand map ${\redptsect{\tr}+\blue{\tilde m}}$ 
we shall now lowerstar along $\red{\langle p,\sect{\tr}\rangle}$
to complete the distributivity pentagon:
\begin{diagram}[w=5.5ex,h=4.5ex,tight]
\redptsect{\tr}+\redsectbluesect{\tr} & \lTo &\cdot& \rTo^{} & \redsectbluesect\tr \\
\dTo<{\redptsect{\tr}+\blue{\tilde m}}  &&&& \dTo>{\blue{\tilde m}} \\
\redptsect{\tr}+\redsect{\tr}&& \rTo_{\red{\langle p,\sect{\tr}\rangle}} && 
\redsect\tr .
\end{diagram}
The resulting right-hand map is just $\blue{\tilde m}$.
Indeed, the lowerstar operation amounts to multiplying along the fibres.
But since in the left-hand summand the map we lowerstar is just the identity of
$\redptsect{\tr}$, no contribution comes from this summand, and in the right-hand
summand we are lowerstarring along the identity map, hence the result is just
$\blue{\tilde m}$, the map that forgets the lower cut.

We have now shown that the set $\redsectbluesect\tr$ appears at the crucial point of 
the diagonal polynomial when constructed from the lower left-hand side.
We proceed to construct it also from the upper right-hand side of the big 
diagram.  The argument is different because of the
non-symmetry of the tensor product.

The pullback square in the upper right-hand corner is, typographically 
transposed:
\begin{diagram}[w=7ex,h=4.5ex,tight]
\redsect{\tr}+\tr&\lTo^{\blue{f+r}}& \sectptsect\tr + \bluesect\tr \\
\dTo<{\red{m}+\tr}&\text{\tiny pb}&\dTo>{\red{\shortsect{m}} +\bluesect\tr}  \\
\tr+\tr&\lTo_{\blue{f+r}}&\blueptsect{\tr}+\bluesect{\tr}.
\end{diagram}
Here there are two new symbols: $\sectptsect\tr$ denotes the set of
trees with a marked cut and a further cut in the tree above the mark.
The map $\red{\shortsect{m}}$ forgets the `short' cut.

We shall now lowerstar the map $\red{\shortsect{m}} +\bluesect\tr$ along
$\blue{f+r}$, to complete the distributivity pentagon
\begin{diagram}[w=5.5ex,h=4.5ex,tight]
\sectptsect\tr + \bluesect\tr  & \lTo &\cdot& \rTo^{} & \redsectbluesect\tr \\
\dTo<{\red{\shortsect{m}} +\bluesect\tr}  &&&& \dTo>{\red{\overline m}} \\
\blueptsect{\tr}+\bluesect{\tr}&& \rTo_{\blue{\langle p,\sect{\tr}\rangle}} && 
\bluesect\tr .
\end{diagram}
The claim is that the right-hand map is just $\red{\overline{m}}$, the map
that forgets the upper cut.  We compute this fibre-wise over an element
$R\subset T $ in $\bluesect\tr$.  The $\blue{\langle p,\sect{\tr}\rangle}$-fibre
has one element for each leaf $e$ of $R$.  For each leaf $e$, the 
$\red{\shortsect{m}}$-fibre consists of the possible cuts in the ideal subtree 
$D_e\subset T$ generated by $e$.  Lowerstarring means multiplying these fibres,
so it amounts to giving a cut in $D_e$ for each leaf $e$ of $R$.  Altogether,
this amounts to giving a total cut in $T$ above the original cut $R \subset T$.
This is once again the set $\redsectbluesect\tr$, but the projection 
$\red{\overline m}:\redsectbluesect\tr\to \bluesect\tr$ this time forgets the 
upper cut.  (Note that there is no contribution from  the right-hand summand,
as it has trivial fibres.)

Finally, note that the square appearing where the two constructions meet,
\begin{diagram}[w=6ex,h=4.5ex,tight]
\redsectbluesect\tr &\rTo^{\red{\overline{m}}}  & \bluesect\tr  \\
\dTo<{\blue{\tilde m}}  &    & \dTo>{\blue{m}}  \\
\redsect\tr  & \rTo_{\red{m}}  & \tr
\end{diagram}
commutes, since it amounts to forgetting first the upper cut and then the lower,
or the other way around.

The remaining part of the big diagram is only about taking pullbacks, and 
presents no difficulties.
\end{blanko}
% 
% \begin{lemma}
%   We have
%   $$
%   M( \tr+m) = \sectsect \tr = M(m+\tr) .
%   $$
%   That is to say that the two constructions `meet' at $\sectsect \tr$ as claimed.
% \end{lemma}

\begin{blanko}{Proof of the unit axiom.}
  The unit is the functor $U:\A^0 \to \A^\tr$ represented by
  $$
  0 \leftarrow 0 \to I \to \tr ,
  $$
  where the last map associates to an edge colour the corresponding trivial tree.
  
  That this is the unit object means that when left-added to the identity functor 
  on $\tr$,
  $$
  0+\tr \lTo 0+\tr \rTo I+\tr \rTo \tr+\tr ,
  $$
  and composing with $M$ yields the identity functor.  And similarly of course
  with right-adding the identity functor.  Checking this amounts to filling the 
  following big diagram:
  \begin{diagram}[w=40pt,h=30pt,tight,scriptlabels]
\tr&&\lTo^\Delta&&\tr&\rTo^\Pi &\tr+I&\rTo^\Sigma&\tr+\tr\\
&\luLig(4,4)&&&\uTo<\Delta&\text{\tiny pb}&\uTo>\Delta&\text{\tiny pb}&\uTo>\Delta\\
\uTo<\Delta&&&&\cdot&\rTo^\Pi&\cdot&\rTo^\Sigma&\ptsect{\tr}+\sect{\tr}\\
&&&&\uTo<\Delta&\text{\tiny pb}&\uTo>\Delta&&\\
\tr&\lTo^\Delta&\cdot&\lTo^\Delta&\tr&\rTo^\Pi&\cdot&\text{\tiny distr}&\dTo>\Pi\\
\dTo<\Pi&\text{\tiny pb}&\dTo<\Pi&\text{\tiny pb}&\dTo<\Pi&\rdLig&\dTo>\Pi&&\\
I+\tr&\lTo_\Delta&\cdot&\lTo_\Delta&\cdot&\rTo_\Pi&\tr&\rTo^\Sigma&\sect{\tr}\\
\dTo<\Sigma&\text{\tiny pb}&\dTo<\Sigma&&\text{\tiny distr}&&\dTo<\Sigma&\rdLig&\dTo>\Sigma\\
\tr+\tr&\lTo_\Delta&\ptsect{\tr}+\sect{\tr}&&\rTo_\Pi&&\sect{\tr}&\rTo_\Sigma&\tr .
\end{diagram}

%     \setlength{\unitlength}{0.7pt}
% 
% \begin{diagram}[w=36pt,h=24pt,tight,objectstyle=\scriptstyle,scriptscriptlabels]
% 0+\tr&&&&&&&&&&&&\\
% &\rdLig(6,4)&&&&&&&&&&&\\
% \uTo<\Delta&&&&&&&&&&&&\\
% &&&&&&&&&&&&\\
% 0+\tr&&\lTo_\Delta&&\cdot&\lTo_\Delta&0+\tr&&&&&&\\
% &&&&&&&\rdLig(2,4)&&&&&\\
% \dTo<\Pi&&\text{pb}&&\dTo<\Pi&\text{pb}&\dTo<\Pi&&&&&&\\
% &&&&&&&&&&&&\\
% I+\tr&&\lTo_\Delta&&\cdot&\lTo_\Delta&\cdot&\rTo_\Pi&\tr&&&&\\
% &&&&&&&&&\rdLig(4,4)&&&\\
% \dTo<\Sigma&&\text{pb}&&\dTo<\Sigma&&\text{distr}&&\dTo<\Sigma&&&&\\
% &&&&&&&&&&&&\\
% \tr+\tr&&\lTo_\Delta&&\ptsect{\tr}+\sect{\tr}&&\rTo_\Pi&&\sect{\tr}&&\rTo_\Sigma&&\tr
% \end{diagram}
% 
% \setlength{\unitlength}{1pt}

As with associativity, we content ourselves to checking that the two 
constructions meet at the distributivity squares.

Starting at the lower left-hand corner with a pullback:

\begin{diagram}[w=7ex,h=4.5ex,tight]
I+\tr&\lTo^{f+r}&V+\sect\tr\\
\dTo<{e+\tr}&\text{\tiny pb}&\dTo>{v + \sect\tr}\\
\tr+\tr&\lTo_{f+r}&\ptsect{\tr}+\sect{\tr}
\end{diagram}
The set $V$ appearing is the set of trees with a marked cut such that 
the mark itself is a leaf (in other words, the upper tree is a trivial
tree).%
\footnote{The author ran out of onomatopoetic notation at this point, and just
chose the letter $V$ at random.}
The map $v: V \to \ptsect\tr$ is the inclusion of these
special pointed cuts into the set of all pointed cuts.

Now lowerstar along $\langle p,\sect{\tr}\rangle$, to complete the
distributivity pentagon:
\begin{diagram}[w=5.5ex,h=4.5ex,tight]
V+\sect\tr  & \lTo &\cdot& \rTo^{} & W&=\tr \phantom{mm} \\
\dTo<{v + \sect\tr}  &&&& \dTo>{w} \\
\ptsect{\tr}+\sect{\tr}&& \rTo_{\langle p,\sect{\tr}\rangle} && 
\sect\tr .
\end{diagram}
This amounts to requiring for each mark of the cut, that the tree above
it is trivial; in other words it is the set $W$ of cuts that only take leaves,
i.e.~the set of trivial root-preserving inclusions, i.e.~the set $\tr$ itself.
(The map $w: W \to \sect\tr$ is of course the inclusion of this particular kind 
of cuts into the set of all cuts.)

The check is similar for the right-hand unit axiom, except that the set in
the corner is the set of trees with the root cut, which again is naturally
identified with $\tr$ itself.
\end{blanko}

\begin{blanko}{Bialgebra of polynomial functors.}
  By the duality of \ref{duality}, the monoidal structure on $\A^\tr$
  induces a comonoidal structure on $\F[\tr]$,
  \begin{diagram}[w=32pt]
  \F  & \lTo^\epsilon   & \F[\tr] & \rTo^\Delta & 
  \F[\tr]\tensor_{\F} \F[\tr]:
  \end{diagram}
  the comultiplication is simply defined by precomposition with $M$ and the
  counit is defined by precomposition with $U$.  In detail, if $F\in\F[\tr]$
  is a polynomial functor in $\tr$-many variables,
  $$
  \tr \lTo E \rTo B \rTo 1,
  $$
  we compute $\Delta(F)$
  as the composite
  \begin{diagram}
  && \ptsect{\tr}+\sect{\tr} & \rTo^{\langle p,\sect{\tr}\rangle}  && \sect{\tr} && && E & \rTo & B &&\\
  & \ldTo^{f+r}&&&&&\rdTo^{m} && \ldTo&&&&\rdTo & \\
\tr+\tr &&&&&&& \tr&&&&&& 1.
\end{diagram}
Note that the comonoidal structure is automatically multiplicative, for formal 
reasons:  multiplication in $\F[\tr]$
is dual to the diagonal functor $\A^\tr\to \A^\tr\times\A^\tr$,
and every monoidal structure on $\A^\tr$ is compatible with this diagonal.
  
It is instructive to calculate $\Delta$ on a multiplicative generator of 
$\F[\tr]$: these are the single-variable polynomials
$$
  \tr \stackrel{\name{T}}\lTo 1 \stackrel{=}\rTo 1 \stackrel{=}\rTo 1,
$$
where the map $\name{T}$ picks out the tree $T\in \tr$.
The composition in this case just amounts to taking fibres, so the result is
  \begin{diagram}[w=6ex,h=4ex,tight]
  && (\ptsect{\tr})_T+(\sect{\tr})_T & \rTo^{\langle 
  p_T,(\sect{\tr})_T\rangle}  && (\sect{\tr})_T && \\
  & \ldTo^{f+r}&&&&&\rdTo & \\
\tr+\tr &&&&&&& 1.
\end{diagram}
Here $(\sect{\tr})_T$ is the set of cuts in the specific tree $T$, 
and the fibre over that is the a two-component set where the left-hand 
component is the set of all pointed cuts on $T$, and the 
right-hand component is again the set of cuts on $T$.
Spelling out what the polynomial does as a functor, we see that it is
exactly the functorial version of the formula from \ref{comultmy}:
we have to sum over the set $(\sect{\tr})_T$ (that's the set of cuts),
then multiply the fibres for each cut $c$; the fibre has an element for each edge 
in the cut (and one extra): for each edge we must apply $f$, which gives the 
branch factor of $P_c$, and multiply all these, getting altogether what is
denoted $P_c$, and for the extra point of the fibre we must apply $r$, which
gives the tree $R_c$ below the cut, which is to be placed in the right-hand 
factor.

The counit associates to each tree $T$ a polynomial in zero variables, 
i.e.~just a finite set.  If the tree is represented as $\name{T}:1 \to \tr$, 
we need to compose
  \begin{diagram}
  && 0 & \rTo  & I && && 1 & \rTo & 1 &&\\
  & \ldTo&&&&\rdTo && \ldTo&&&&\rdTo & \\
0 &&&&&& \tr&&&&&& 1 .
\end{diagram}
The result is
  \begin{diagram}
  && 0 & \rTo  && I \times_\tr 1 && \\
  & \ldTo&&&&&\rdTo & \\
0 &&&&&&& 1 ,
\end{diagram}
i.e.~the set $I\times_\tr 1$, which is the 
singleton set $1$
if $T$ belongs to $I$ and the empty set $0$ otherwise. Again it is
clear that this is the functorial analogue of the specification in 
Section~\ref{Sec:bialg}.
\end{blanko}

\section{Beyond $P$-trees}
%%%%%%%%%%%%%%%%%%%%%%%%%%%%%%%%%%%%%%%%%%%%%%%%%%
\label{Sec:beyond}

\begin{blanko}{Trees (irrespective of any $P$).}
  We have worked with $P$-trees for two reasons: the first is to allow the
  connection with the free-monad construction.  The second is the fact that
  $P$-trees (in the set-based setting we are working in) are rigid.  This
  guarantees that taking sets of isomorphism classes is well-behaved.

  Abstract trees, in the sense of \ref{polytree-def}, are not $P$-trees for any
  polynomial endofunctor over $\Set$.  The `terminal endofunctor' in the 
  category of polynomial endofunctors and their cartesian morphisms \eqref{equ:cartmorphism}
  which should make the abstract trees $P$-trees does not exist.  Also, trees
  clearly may have non-trivial automorphisms, so taking sets of isomorphism
  classes may be prone to errors, and care is needed.
\end{blanko}
  
\begin{blanko}{Reinterpretation of the sets and maps in Section~\ref{Sec:monoidal}.}
  In spite of the above remarks, it is nevertheless possible to mimic all the
  constructions of the previous section, if just the involved sets and set maps
  are defined a little bit differently.
  
  The starting point is the same: we now let $\tr$ denote the set of
  isomorphism classes of abstract trees.  Hence the basis for the bialgebra is
  the same as used already in Section~\ref{Sec:bialg}.
  But while for $P$-trees, the set $\sect\tr$ can be described as isomorphism
  classes of root-preserving inclusions, this description does not work for
  abstract trees.  Indeed, it is crucial that the fibres of the projection $m:\sect\tr\to\tr$
  have the correct cardinality: for example, the abstract tree
  \begin{center}\begin{texdraw}
    \move (0 -12) \lvec (0 0)\onedot \lvec (-8 10)\onedot \lvec (-12 22) 
    \move (-8 10) \lvec (-4 22)
    \move (0 0) \lvec (8 10)\onedot \lvec (12 22) \move (8 10) \lvec 
    (4 22)
  \end{texdraw}\end{center}
  certainly has five different cuts, but two of them are isomorphic as abstract 
  cuts.  So instead of defining $\sect\tr$ as the set of isomorphism classes
  of cuts, it is necessary first to choose a representative for each iso-class
  in $\tr$, and then for each such representative $T$ define $m^{-1}$ to be the set of
  cuts of that specific tree $T$, and finally let $\sect\tr$ be the disjoint union
  of all those fibres.   Similar care is needed to define $\tr'$ and 
  $\ptsect\tr$, and the other symbols involved.  The fact that all the symbols
  are defined fibrewise over $\tr$, instead of being defined abstractly in
  terms of isomorphism classes, guarantees that all the pullback and lowerstar
  operations (which are fibrewise operations)  yield the expected results,
  and one can check that the proof of associativity and the unit axiom work
  equally well for abstract trees, with these provisos.
\end{blanko}

\begin{blanko}{Groupoids instead of sets.}
  Clearly the remedies just explained are rather ad hoc.  A better and more 
  conceptual way to account for abstract trees and $P$-trees on equal footing
  consists in upgrading the theory from sets to groupoids.  This can be seen
  as one further step of categorification.
  
  First of all it is necessary (and possible) to upgrade the theory of
  polynomial functors from sets to groupoids: this means that the representing
  diagrams for polynomial functors 
  $$
  I \lTo E \rTo B \rTo J
  $$
  should then be (suitably homotopically finite) groupoids, and the functors
  themselves go between slices of the category of groupoids (slices being taken
  in the homotopical sense).  With the appropriate adjustments, everything works
  as for sets: pullbacks and fibres have to be homotopy pullbacks and homotopy
  fibres, if the maps are not groupoid fibrations; sums are replaced by slightly
  fancier colimits, and everything is up to equivalence of groupoids instead of
  isomorphism of sets.  This theory correctly incorporates all questions of
  symmetries of objects \cite{Kock:MFPS28}.
  
  In this setting, abstract trees are $E$-trees for the
  `exponential functor' 
  $$
  E(X)= \sum_{n\in\N} X^n/\Aut[n],
  $$ represented by the
  groupoid diagram
  $$
  1 \lTo \B' \rTo \B \rTo 1
  $$
  where $\B$ is the groupoid of finite sets and bijections, and $\B'$ the
  groupoid of finite pointed sets and basepoint-preserving bijections, and where
  the quotient is the weak groupoid quotient (sewing in paths instead of 
  collapsing).   See Baez--Dolan~\cite{Baez-Dolan:finset-feynman} for further
  discussion of such constructions, and the exponential functor in particular.
  
  In the groupoid setting, one works with the groupoid of trees, instead of with
  its set of isomorphism classes.  The practical benefit, in the presence of
  automorphisms of objects, is the same as working with moduli stacks rather
  than coarse moduli spaces in algebraic geometry.

  The choice in this paper to work with sets is first of all that this theory
  is already available.  Second, it is preferable to expose the essential ideas
  and constructions in the transparent rigid set-up, so as not to burden the
  arguments with homotopy theory.  However, to ease the eventual upgrade to
  groupoids, the exposition is careful to perform the arguments in a clean,
  functorial language, favouring natural isomorphisms over element-based
  constructions.
  
  In fact, the theory of polynomial functors is currently being worked out for
  infinity-groupoids, in joint work with David Gepner.  One important insight is
  that $\infty$-groupoids play just the same role for $\infty$-categories as
  sets play for categories.  Seeing how the theory goes for sets is therefore an
  important first step for the more general case of groupoids and
  $\infty$-groupoids.
\end{blanko}

% \nocite{Carboni-Lack-Walters}
% \nocite{Kock:NotesOnPolynomialFunctors}
\nocite{Connes-Kreimer:9808042}
\nocite{Connes-Kreimer:9912092}

\hyphenation{mathe-matisk}

% \bibliographystyle{scplain}
% \bibliography{joachims}

\label{lastpage}
\end{document}